\documentclass[11pt]{article}

\pdfoutput=1

\usepackage{amssymb}
\usepackage{amsmath}
\usepackage[affil-it]{authblk}
\usepackage[backend=biber,citestyle=numeric-comp,bibstyle=ieee,dashed=false]{biblatex}
\usepackage{caption}
\usepackage{enumitem}
\usepackage{fancyhdr}
\usepackage{float}
\usepackage[bottom,hang]{footmisc}
\usepackage[a4paper,left=2.5cm,right=2.5cm,top=2.5cm,bottom=2cm]{geometry}
\usepackage{graphicx}
\usepackage{microtype}
\usepackage{newtxtext}
\usepackage{placeins}
\usepackage{scalerel}
\usepackage{tabularx}
\usepackage[colorlinks=true,allcolors=blue]{hyperref}

\def\equationautorefname~#1\null{(#1)\null}
\AtEveryBibitem{\ifentrytype{article}{\clearfield{issn}}{}}
\DefineBibliographyExtras{english}{}
\DefineBibliographyStrings{english}{phdthesis = {Ph.D.~thesis}}
\DeclareCiteCommand{\cite}
	{\usebibmacro{cite:init}%
	 \usebibmacro{prenote}}
	{\usebibmacro{citeindex}%
 		\ifnumequal{\value{citecount}}{1}{\printtext[bibhyperref]{[}}{}%
		\usebibmacro{cite:comp}%
	}
	{}
	{\usebibmacro{cite:dump}%
	 \ifnumequal{\value{citecount}}{\value{citetotal}}{\printtext[bibhyperref]{]}}{}%
	 \usebibmacro{postnote}%
	}

\DisableLigatures[f]{encoding = *, family = *}
\stdpunctuation
\let\Oldsection\section
\let\Oldsubsection\subsection
\let\Oldsubsubsection\subsubsection
\renewcommand{\section}{\FloatBarrier\Oldsection}
\renewcommand{\subsection}{\FloatBarrier\Oldsubsection}
\renewcommand{\subsubsection}{\FloatBarrier\Oldsubsubsection}

\newsavebox{\foobox}
\newcommand{\slantbox}[2][0]{\mbox{%
        \sbox{\foobox}{#2}%
        \hskip\wd\foobox
        \pdfsave
        \pdfsetmatrix{1 0 #1 1}%
        \llap{\usebox{\foobox}}%
        \pdfrestore
}}
\newcommand\unslant[2][-.3]{\slantbox[#1]{$#2$}}

\newcommand{\titleOfPaper}{Vibration Reduction by Stiffness Modulation\\--\\a Theoretical Study}
\author[a,$\ast$]{Alexander Nowak}
\author[a]{L. Flavio Campanile}
\author[a]{Alexander Hasse}
\newcommand{\mailA}{alexander.nowak@mb.tu-chemnitz.de}

\affil[a]{Chemnitz University of Technology\par
		  Professorship Machine Elements and Product Development, 
		  Reichenhainer Str. 70, 09126 Chemnitz, Germany\vspace{3pt}}
\affil[$\ast$]{Corresponding author. E-mail address: \href{mailto:\mailA}{\mailA}}

\title{\textbf{\titleOfPaper \vspace{1em}
	   }}
\date{}
\begin{document}
\maketitle
\vspace{-20pt}
\hrule

\begingroup
\renewcommand{\thefootnote}{}
\footnotetext[1]{
	\textcopyright\ 2021. This manuscript version is made available under the \href{http://creativecommons.org/licenses/by-nc-nd/4.0/}{CC-BY-NC-ND 4.0 license}.\\
	Published journal article: \textit{Journal of Sound and Vibration}, vol. 501, pp. 116040, 2021. \textsc{doi}: \href{https://doi.org/10.1016/j.jsv.2021.116040}{10.1016/j.jsv.2021.116040}.\\
	The article can be accessed at \url{https://www.sciencedirect.com/science/article/pii/S0022460X21001127}.
	}
\setcounter{footnote}{0}
\renewcommand{\thefootnote}{\value{footnote}}
\endgroup

\section*{\small \centering Abstract}
\label{sec*:abstract}
\begingroup
\small
\leftskip1cm
\rightskip\leftskip
Semi-active vibration reduction techniques are defined as techniques in which controlled actions do not operate directly on the system's degrees of freedom (as in the case of active vibration control) but on the system's parameters, i.e., mass, damping, or stiffness.

Cyclic variations in the stiffness of a structural system have been addressed in several previous studies as an effective semi-active vibration reduction method.
The proposed applications of this idea, denoted here as \textit{stiffness modulation}, range from stepwise stiffness variations on a simple spring-mass system to continuous stiffness changes on rotor blades under aerodynamic loads.

Semi-active systems are generally claimed to be energetically passive.
However, changes in stiffness directly affect the elastic potential energy of the system and require external work under given conditions.
In most cases, such injection or extraction of energy (performed by the device in charge of the stiffness variation and denoted as the \textit{pseudo-active} effect) usually coexists with the \textit{semi-active} effect, which operates by redistributing the potential energy within the system in such a way that it can be dissipated more efficiently.

This work focuses on the discrimination between these two effects, which is absent in previous literature.
A first study on their dependence on the process parameters of stiffness modulation is presented here, with emphasis on the spatial distribution of the stiffness changes.
It is shown that localized changes tend to result in a larger semi-active share of vibration attenuation, whereas a spatially homogeneous stiffness modulation only generates a pseudo-active effect.

\begin{enumerate}[label=\textbf{Keywords:},leftmargin=*]
\leftskip1cm
\rightskip\leftskip
\item {
	Stiffness Modulation, Variable Stiffness, Switched Stiffness, Semi-Active Vibration Reduction, Vibration Control.
	}
\end{enumerate}
\endgroup
\vspace{8pt}
\hrule

\section{Introduction}
\label{sec:1}
\subsection{State of the art}
\label{subsec:1.1}
Structures and machines are, in general, capable of experiencing vibrations owing to the mutual interaction between elastic and inertial forces.
In cases where vibrations are unwanted and the inherent damping of a given system is not sufficient to keep the decay time and resonance amplitude below acceptable limits, vibration suppression techniques come into play to modify the system's dynamics accordingly.

A multiplicity of approaches for vibration suppression are known from the literature, which can be classified in different ways.
The most commonly used classification subdivides vibration suppression approaches into the categories \textit{passive}, \textit{active}, and \textit{semi-active}.
Passive methods are those that do not require actuators or other controlled devices.
Active methods interact with the vibrating system by controlled actuators which directly operate on the system's degrees of freedom (DoF).
Semi-active methods involve controlled devices that act on the parameters of the system \cite{Nitzsche2012}.
This paper focuses on semi-active methods that modify the system's stiffness.

Semi-active methods that employ stiffness changes to influence vibrations can be further classified according to two distinct principles.
The first principle (\textit{adaptive stiffness}) uses variable stiffness to modify, in a quasi-static manner, the system's transfer function, thus improving its response to external disturbances.
The second principle uses cyclical stiffness changes related to the vibration signal.
Following the wording of Anusonti-Inthra and Gandhi \cite{Anusonti-Inthra2000}, it involves \textit{modulating} the stiffness of the system in real time and is therefore referred to here as \textit{stiffness modulation}.

Adaptive stiffness systems can modify their resonance frequencies \cite{Greiner-Petter2014} which can thus be moved to values with a low spectral amplitude of the excitation.
Stiffness adaptation can also be used to change the resonance frequency of tuned absorbers or dampers to match the excitation frequency of the main system \cite{Franchek1996}.
Electromechanical stiffness adjustment of absorbers by means of shunted piezoceramic elements is covered in \cite{Davis2000}.
In \cite{Lin2015}, a tuned mass damper is supplemented with a resettable variable stiffness device which adapts the system's stiffness and damping.
Other concepts for adaptive devices include the use of magnets with adjustable mutual distance \cite{Sayyad2014}, conical springs \cite{Suzuki2013}, actuator-controlled variation of the angle between springs \cite{Nagarajaiah2005}, magnetorheological elastomers \cite{Komatsuzaki2015}, and variable-length pendula \cite{Pasala2014,Wang2020}.

If the stiffness is changed in an actively controlled and cyclical manner in periods comparable to the period of the vibration under consideration (stiffness modulation), a significant reduction in the vibration amplitude can be achieved.
Main contributions on this topic focus on building seismic isolation \cite{Xinghua2000,Nagarajaiah2006,Tan2004,Liu2006} and helicopter vibration suppression \cite{Anusonti-Inthra2000,Yong2004}.
Studies on shock isolation are presented in \cite{Ledezma-Ramirez2011,Ledezma-Ramirez2012,Ledezma-Ramirez2014}.
General purpose and basic studies can be found in \cite{Liu2008,Corr2001,Onoda1992} and \cite{Ramaratnam2006}.

\begin{figure}[hbt]
	\centering
	\includegraphics{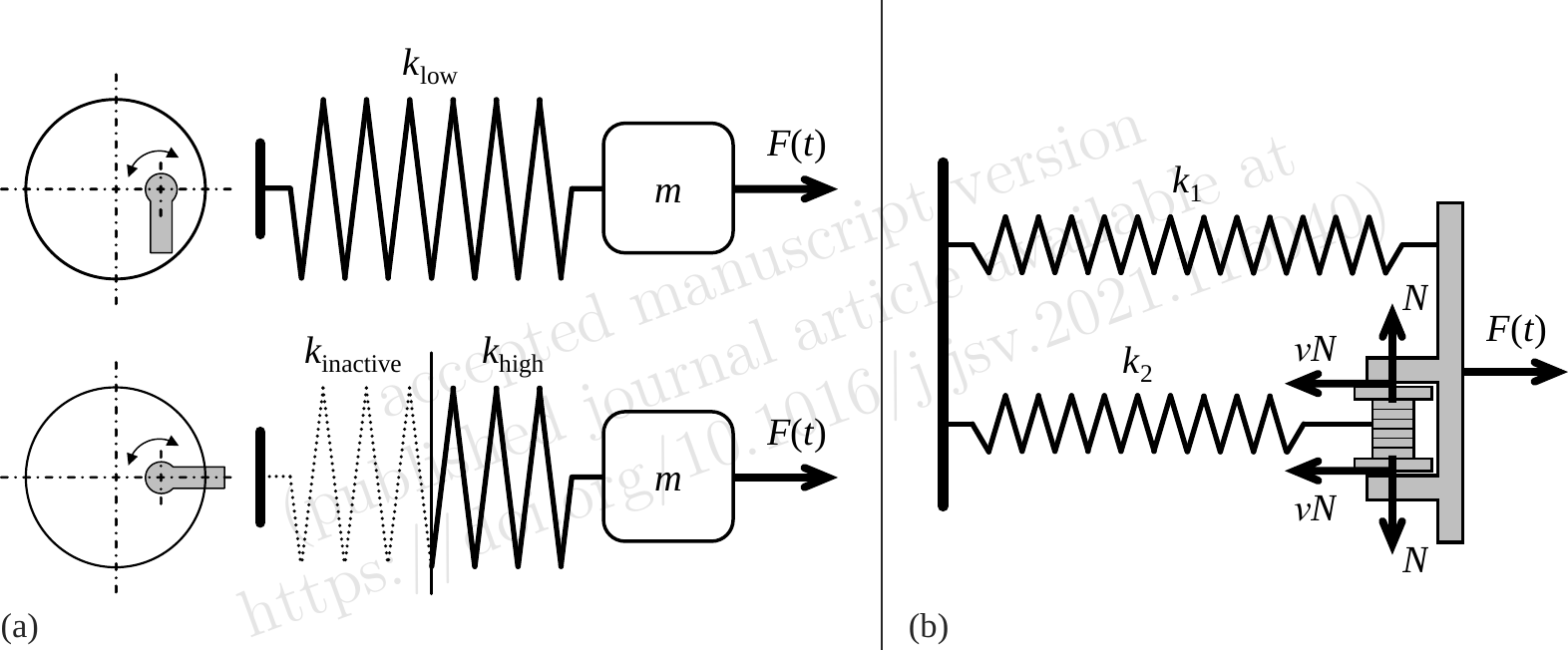}
	\caption{(a) Spring with adaptable active length according to \cite{Ramaratnam2006} and (b) ``Smart Spring'' according to \cite{Yong2004}}
	\label{fig:01}
\end{figure}

In the context of stiffness modulation, the logic that rules the stiffness changes as a function of vibration signals is called the \textit{control logic}.
A first control logic concept (switched stiffness), adopted among others in \cite{Leitmann1994,Corr2001,Ramaratnam2006}, uses two possible values for the stiffness.
At the extremal points of the vibration (maximum absolute value of the displacement and zero velocity), the stiffness is switched from the higher to the lower value, whereas at zero crossings the higher stiffness value is restored.
In \cite{Corr2001}, the stiffness change is realized by piezoelectric actuators, which can be switched from an open-circuit state (high stiffness) to a short-circuit state (low stiffness).
In \cite{Ramaratnam2006}, the stiffness of a spring-mass oscillator is varied using a motor-controlled arm, which changes the effective spring length (see \autoref{fig:01}a).
Depending on the position of the arm, some coils can be blocked and rendered inactive, which results in increased stiffness.
A similar control logic is adopted in \cite{Liu2006} and \cite{Onoda1992}.
The first contribution deals with magnetorheological dampers, in which the damping is varied along with the stiffness.
The second study investigates a tension-controlled string.
The aforementioned studies on shock isolation \cite{Ledezma-Ramirez2011,Ledezma-Ramirez2012,Ledezma-Ramirez2014} also adopt the switched-stiffness control logic.

In \cite{Xinghua2000}, the primary structure is connected by a hydraulic cylinder to an elastic brace.
When the cylinder control valve is closed, the stiffness of the brace acts on the primary structure;
by opening the valve, the value of the stiffness connected to the primary system drops to zero.
Other than in the switched-stiffness case, in which the low-stiffness state is maintained until the zero crossings, in \cite{Xinghua2000}, the stiffness is reduced for a short time about the extremal points of the vibration.
The hardware concept presented in \cite{Tan2004} is virtually identical, except for the presence of additional passive dampers.
The system is driven by an on-off controller.

The ``Smart Spring'' presented in \cite{Yong2004} (see \autoref{fig:01}b) works according to a principle similar to that in \cite{Xinghua2000} and \cite{Tan2004}, since it adds a spring in a parallel configuration to the primary system.
The coupling is realized by friction forces controlled by a piezoceramic actuator.
Strictly speaking, also this solution is limited to two possible stiffness values (primary structure with or without the additional spring).
Intermediate values of the normal force produced by the actuator generate stick-slip effects with piecewise linear force-displacement curves.
The authors assign to these states intermediate stiffness values, which are linearly interpolated with respect to the electric voltage applied to the piezoceramic actuator.
The voltage is varied according to an adaptive control law.

The device presented in \cite{Ledezma-Ramirez2011,Ledezma-Ramirez2012,Ledezma-Ramirez2014} is similar to the Smart Spring, but without a friction interface.
In \cite{Ledezma-Ramirez2014}, the main dissipation mechanism is explained as the result of the inelastic impact between the masses associated with the two degrees of freedom at the time of coupling (switching from low to high stiffness).
In \cite{Ledezma-Ramirez2011}, the presence of a viscous damper is considered instead.
A practical realization of the stiffness variation by means of magnets is addressed in \cite{Ledezma-Ramirez2012}.
A similar setup is used in \cite{Tran2017} to study the effect of time delay on variable stiffness control.

In \cite{Liu2008}, the serial arrangement of a spring and a Voigt element with a controllable damper is used.
In the experiment, a magnetorheological damper is used for this purpose.
For low values of the damping coefficient, a low-stiffness condition is reached, in which the stiffness of the spring and of the Voigt element work in series.
If the controllable damper is set to high damping, then the Voigt element tends to be rigid and the stiffness of the system approaches the value of the spring (high-stiffness condition).
An on-off control scheme is used to select between the two configurations.

In \cite{Anusonti-Inthra2000}, an extensive study on the effect of cyclic variations of the root stiffness on the vibration of rotor blades is presented.
No concept is provided for the hardware interface expected to produce stiffness variations, which are considered as given and processed numerically.
The bending stiffness along the two axes and the torsional stiffness are varied independently and in a sinusoidal fashion.
Different frequencies (as multiples of the rotational frequency) and phase lags are investigated.

In \cite{Nagarajaiah2006}, a device similar to \cite{Nagarajaiah2005} is used for base isolation employing stiffness modulation.
The system with four springs in a rhombus configuration is used to influence the stiffness of the system by changing the shape of the rhombus.
The springs are stabilized against buckling by telescopic tubes, which also add friction forces to the system.
Different time signals (sinusoidal, square, triangle, and random) are used to drive the stiffness-variable system.

A final mention concerns a concept in which the stiffness variation device is not actively controlled but mechanically coupled with the vibrating system.
In \cite{Anubi2013}, a secondary vibrating system consisting of a ``control mass,'' a restoring spring, and a damper, is connected to the stiffness variation device of the main system.
In a later work, the possibility of influencing the motion of the control mass by a magnetorheological damper is investigated \cite{Anubi2015}.
In this last option, the parameters of a secondary vibrating system are modified, which in turn influences through its motion the parameters of the primary system.

\subsection{Semi-active vs. pseudo-active}
\label{subsec:1.2}
As previously mentioned, vibration control techniques are defined as semi-active if they act on the system's parameters.
Although a large number of contributions dealing with semi-active techniques neither provide a definition of the ``semi-active'' attribute nor a corresponding citation, the proposed solutions meet the classification criterion mentioned above.

However, there is a relatively large number of contributions in which the ``semi-active'' attribute is related to energy considerations, which can lead to ambiguities.
In some cases, semi-active measures are said to require less energy to operate \cite{Tan2004,Liu2005,Liu2006,Liu2008,Clark2000}.
According to Xingua \cite{Xinghua2000}, semi-active devices have extremely low power requirements.
Suzuki and Abe \cite{Suzuki2013} went even further and asserted that a semi-active control system does not need a power source to suppress vibration, while Fisco and Adeli \cite{Fisco2011} claimed that it ``is normally operated by battery.''
More substantial statements deal with the capability to perform work on the system.
According to Anubi and Crane \cite{Anubi2015}, semi-active devices are either dissipative or conservative.
Similarly, Liu \cite{Liu2004} asserted that they ``can only dissipate energy'' and ``cannot put energy into the system.''
The second assertion can also be found in other studies \cite{Leitmann1993,Spencer1997,Spencer2003,Pasala2013}.

Following this more restrictive definition, a vibration control device acting on the system's parameters but capable of performing work on the system cannot be classified as semi-active, and would have to be labeled as an active control device.
A variable-stiffness device, however, mostly requires the ability to perform work on the system, since a change in the stiffness of a vibrating system modifies its elastic potential energy.
Indeed, while several authors agree with the fact that semi-active measures, owing to the above-mentioned property of being unable to perform work on the system, cannot destabilize it \cite{Spencer2003,Pasala2013,Preumont2011}, in the special case of variable-stiffness devices it is recognized that instabilities can occur \cite{Winthrop2004,Corless1997}.

Within this more restrictive definition, it remains unclear whether the vibration control device must not perform positive work on the system to be classified as semi-active, or if any kind of energy exchange (positive or negative) disqualifies it from this category.
As vibration attenuation implies a reduction in the energy level, the case of negative work should be the rule and needs to be considered.
Performing negative work implies loading the system with controlled forces, which requires, in turn, a proper actuator and energy supply.
Therefore, it can be expected that the resulting cost, effort, and complexity primarily depend on the managed energy flow, and only secondarily on its sign.

Based on these considerations, we refer to the most restrictive choice and define as (purely) \textit{semi-active} a vibration suppression effect that acts on the system's parameters (especially on the stiffness), and does not add or subtract energy through the stiffness variation device.
If the stiffness variation device performs (positive or negative) work on the system, this effect is referred to here as \textit{pseudo-active}.
Because no energy can flow through the stiffness variation device with the semi-active effect, this can only involve rearranging the vibration energy in the system so that it is better dissipated by internal damping.
The semi-active effect is therefore quantified after subtracting the energy reduction due to internal damping, which would be present without stiffness modulation (\textit{passive} effect).

\subsection{Motivation of this study}
\label{subsec:1.3}
Analysis of published work leads to following considerations:
\vspace{-0.5\topsep}
\begin{itemize}
\setlength{\itemsep}{-1pt}
	\item
A large part of the considered contributions on stiffness modulation deals with isolators \cite{Nagarajaiah2006,Liu2006,Ledezma-Ramirez2011,Ledezma-Ramirez2012,Ledezma-Ramirez2014,Liu2005,Liu2004} or interfaces acting on the boundary conditions of the vibrating system \cite{Anusonti-Inthra2000,Yong2004}.
Hence, all these approaches do not consider the variation of the internal stiffness of the vibrating system.
A second large group of studies involves principle investigations on single-mass oscillators;
practical applications or the extension to more complex structures is not discussed \cite{Leitmann1994,Corr2001,Ramaratnam2006,Tran2017}.
Solutions involving stiffness changes of a structural system of practical relevance and complexity are treated in \cite{Xinghua2000,Tan2004}.
	\item
All investigated practical implementations of variable stiffness involve lumped elements or forces (springs, dampers, hydraulic cylinders, piezoceramic stack actuators, magnets), which are primarily devised to influence a single degree of freedom of the system.
	\item
Many solutions \cite{Greiner-Petter2014,Nagarajaiah2006,Tan2004,Liu2006,Yong2004,Corr2001,Anubi2013,Anubi2015} employ additional dissipative elements, such as dampers, friction interfaces, or resistors.
This makes it difficult to judge the potential of mere stiffness modulation.
\end{itemize}
A need for dedicated investigations can be identified with a focus on the following points:
\vspace{-0.5\topsep}
\begin{itemize}
\setlength{\itemsep}{-1pt}
	\item
Stiffness modulation involving stiffness changes of the main load-carrying, vibration-prone system (and not of an external interface or isolator);
	\item
Stiffness modulation involving distributed stiffness changes;
	\item
Stiffness modulation that relies on the inherent damping properties of the vibrating system as dissipation mechanism and does not require additional damping elements.
\end{itemize}
With this long-term scenario in mind, the present work focuses on the distinction between the above defined semi-active and pseudo-active effects, and on identifying criteria to maximize the semi-active part of the vibration reduction achieved by stiffness modulation.
A large semi-active part leads to a low power flow through the stiffness variation devices.
These can be built, in turn, in a more compact way, and therefore be more easily integrated and/or distributed in the main load-carrying system.

This paper presents a theoretical study on stiffness modulation using switching stiffness control logic.
The stiffness changes are numerically imposed on the considered lumped-parameter systems without modeling a physical stiffness variation device.
The analysis focuses on the above introduced energy discrimination, distinguishing between pseudo-active, (purely) semi-active, and passive effects.
As an essential parameter of this analysis, the effect of the localization of stiffness changes is investigated.
After an illustrative part based on spatial coordinates, the analysis proceeds in the modal space, and free as well as forced oscillations are studied.

\section{Investigation on spring-mass oscillators}
\label{sec:2}
\subsection{Switching stiffness control logic}
\label{subsec:2.1}
In \cite{Leitmann1994}, the already introduced switching stiffness control logic was identified, for a single-DoF system, as the control logic that maximizes the energy extracted from the system.
It can be formulated as follows:
\begin{equation}
	k = \begin{cases}
			k_h,	& \text{if}~~y(t)\dot{y}(t) \geq 0 \\
			k_l,	& \text{if}~~y(t)\dot{y}(t) < 0
		\end{cases}
	\label{eq:01}
\end{equation}
where $ k_h $ and $ k_l $ are the high and low values of the system's stiffness, respectively, and $ y(t) $ is the vibration signal.

To include systems with more than one degree of freedom, we specify the stiffness for any single stiffness parameter~$ j $ via a corresponding factor~$ \gamma_j $
\begin{equation}
	k_j = \begin{cases}
			\gamma_j k_{0j},	& \text{if}~~c(t)\dot{c}(t) \geq 0 \\
			k_{0j},				& \text{if}~~c(t)\dot{c}(t) < 0
		\end{cases}
	\label{eq:02}
\end{equation}
which scales the stiffness between low and high stiffness.
The value $ c(t) $ is a properly chosen \textit{observation variable} of the vibrating system.

To obtain a criterion for real-time control, it can be observed that the upper part of condition~\autoref{eq:02} is verified when the absolute value of the observation variable increases, while the lower part is verified when it decreases.
It can then be written in terms of two consecutive points in time $ t_i $ and $ t_{i+1} $ as
\begin{equation}
	k_j = \begin{cases}
			\gamma_j k_{0j},	& \text{if}~~|c(t_{i+1})| \geq |c(t_{i})| \\
			k_{0j},				& \text{if}~~|c(t_{i+1})| < |c(t_{i})| 
		\end{cases}
	\label{eq:03}
\end{equation}

In the case in which the observation variable oscillates around zero with only one extremal point between two zeroes, the stiffness is increased directly after each zero crossing and is decreased directly after each extremal point.
In the following, the stiffness decrease and increase points are labeled as $ e_i $ and $ z_i $, respectively.

All following considerations are based on the switched-stiffness control logic presented here.

\subsection{Coupleable two-DoF system}
\label{subsec:2.2}
\begin{figure}[hbt]
	\centering
	\includegraphics{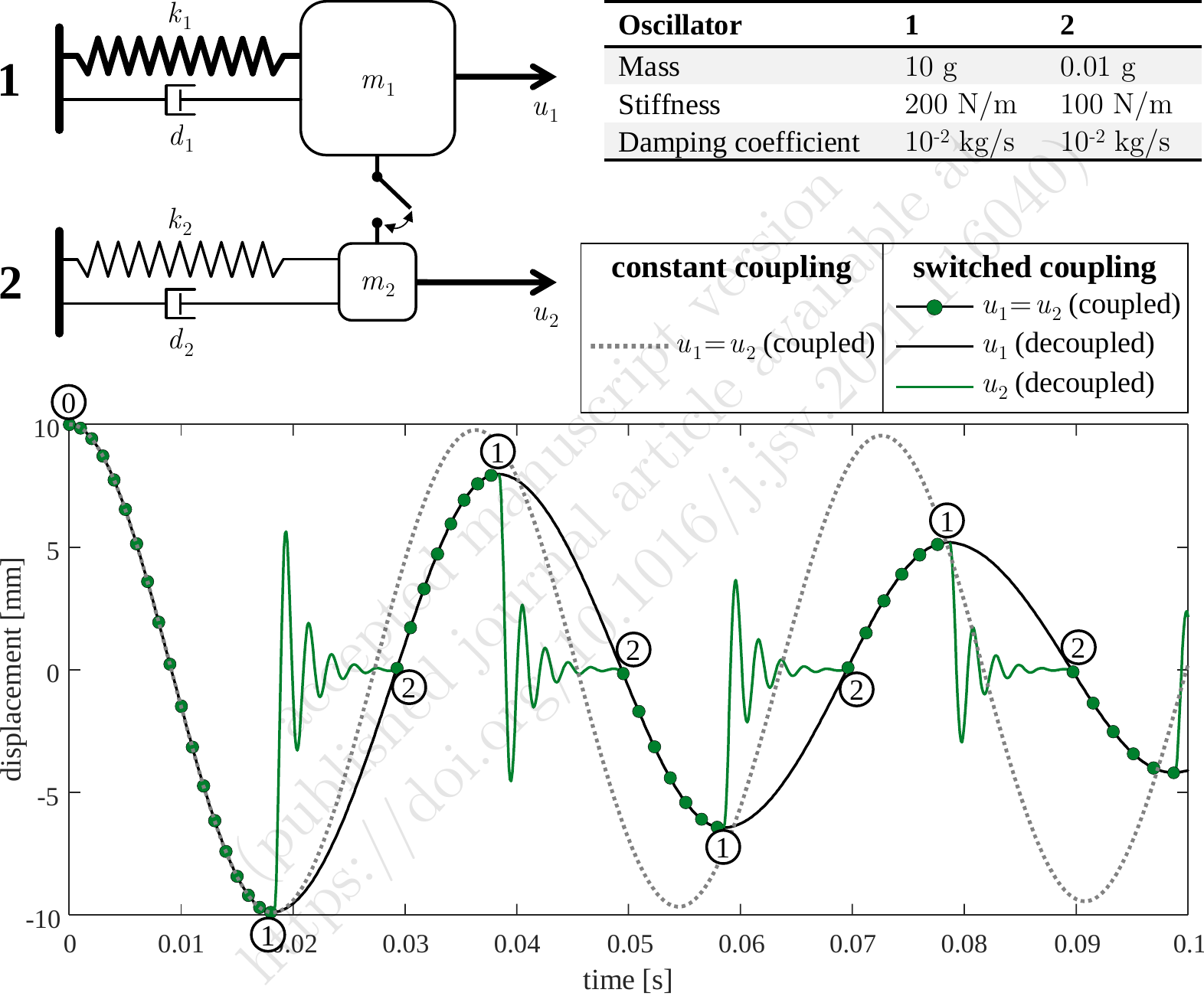}
	\caption{Switched stiffness approach on a system of two coupleable oscillators}
	\label{fig:02}
\end{figure}

In many systems studied in the literature (e.g., \cite{Yong2004} -- see \autoref{fig:01}b -- and \cite{Xinghua2000,Tan2004,Ledezma-Ramirez2011,Ledezma-Ramirez2012,Ledezma-Ramirez2014}) the stiffness change is achieved by coupling and decoupling the vibrating system with an additional stiffness in parallel.
As a first illustrative example (see \autoref{fig:02}), such an option is studied here, where the secondary system also possesses its own mass.
However, the mass of the secondary system is relatively small;
therefore its effect on the dynamics of the primary system in the coupled state is of limited relevance.
The two oscillators have the same damping coefficient ($d_1=d_2$), and the coupling connection is assumed to be massless and ideally stiff.
A case is considered in which the masses $m_1$ and $m_2$ initially oscillate in a coupled configuration (see \autoref{fig:02}, Point~0 as the starting point for the combined oscillation, defined by $u_1=u_2$).
After the first zero-crossing, the switching control logic \autoref{eq:01} is applied.
At the point of maximum deflection of the two coupled springs, the system is decoupled (see \autoref{fig:02}, Point~1).
As a result, the small spring begins to oscillate at a higher frequency with a rapidly decreasing amplitude (see $u_2$ curve in decoupled configuration) such that this oscillation can be assumed to be extinguished before the springs are coupled again at the zero crossing (see \autoref{fig:02}, Point~2).
Through the coupling, the remaining kinetic energy of the system, stored in the first mass, is redistributed between the two masses.
In the subsequent cycles, the process is repeated, resulting in the progressive decay of the amplitude of oscillator~1.

Recalling \autoref{eq:01}, $ k $ corresponds in the present case to the stiffness acting on the first mass, i.e., $ k_l $ is the stiffness of the first spring and $ k_h $ is the stiffness of the parallel connection $ k_1+k_2 $.
The observation variable $ c $ is given by the first mass displacement $ u_1 $.

The example in \autoref{fig:02} shows that the vibration reduction of the overall system is significantly faster with cyclic coupling and decoupling of the springs (see graph in \autoref{fig:02}, \textit{switched coupling}) than with permanently coupled springs (see graph in \autoref{fig:02}, \textit{constant coupling}).
Damping has a much stronger effect on the high-frequency oscillator 2 than on the coupled state.
The switched stiffness approach exploits this effect by repeatedly shifting energy to oscillator~2 to dissipate it more effectively.

\subsection{Single-DoF system}
\label{subsec:2.3}
Let us now consider the case of an undamped single-DoF system with a constant mass $ m = m_1 $ and variable stiffness according to \autoref{eq:01} with
\begin{equation}
	k_l = k_1, \quad
	k_h = k_1 + k_2
	\label{eq:04}
\end{equation}

\begin{figure}[hbt]
	\centering
	\includegraphics{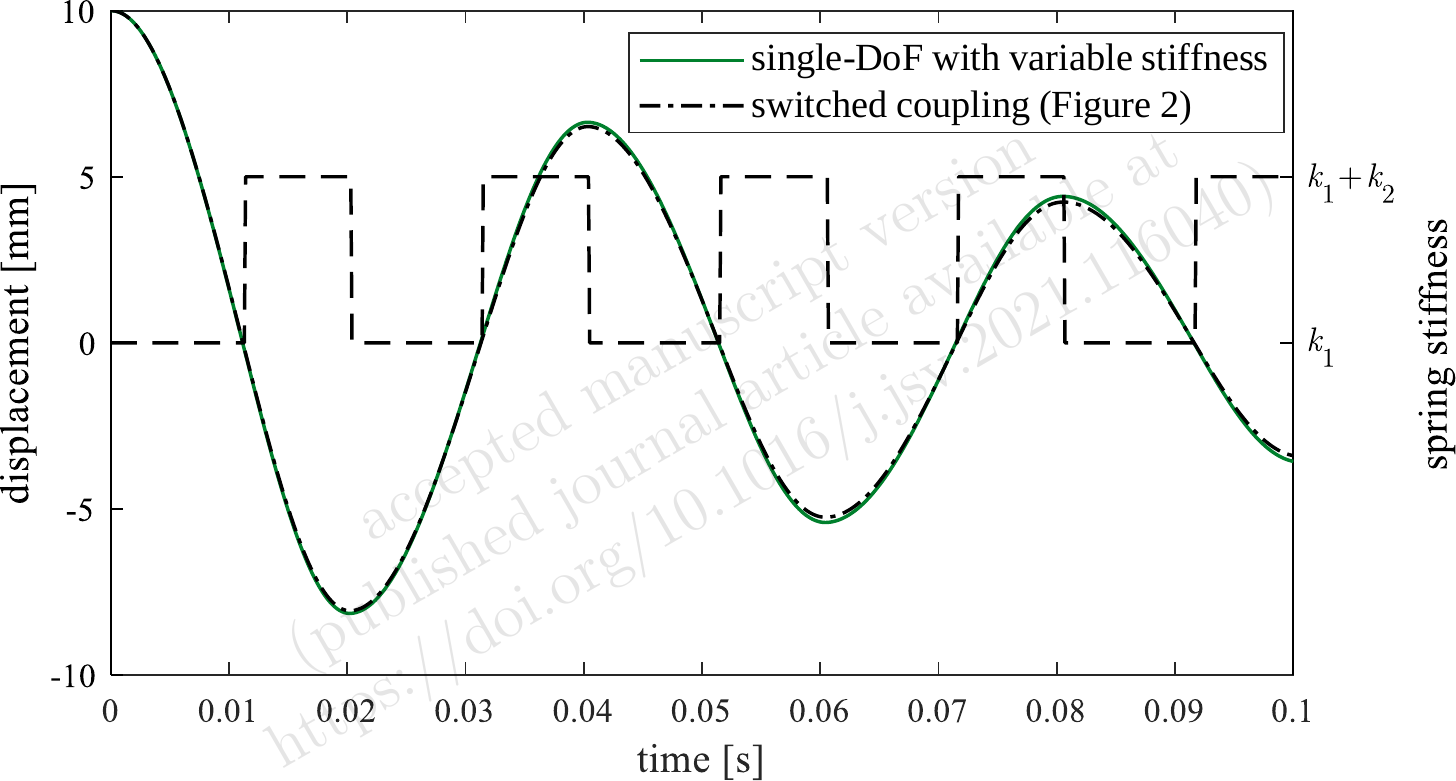}
	\caption{Single-DoF oscillator compared with the coupleable two-DoF system in \autoref{fig:02}}
	\label{fig:03}
\end{figure}

By changing the stiffness along a square wave as shown in \autoref{fig:03} (dashed curve), nearly the same decay curve (see \autoref{fig:03}, green curve) as in the first example (see \autoref{fig:03}, black curve) can be obtained.
There is only a slight difference owing to the invariable mass $ m $ and the absence of damping.
Methodically, however, the two examples differ essentially from each other:
the energy that was transferred to another passive subsystem (and then dissipated) in the first case is now extracted from the system at the points of stiffness reduction.
The stiffness variation device performs negative work on the system.
This negative work (called \textit{extracted energy} in the following) is equal to the variation in the potential energy
\begin{equation}
	\Delta U = \dfrac{1}{2} \Delta k u^2
	\label{eq:05}
\end{equation}
which results from the stiffness change $ \Delta k $ at displacement $ u $.

Hence, in the first case, the vibration-reducing effect is purely semi-active, as defined in \autoref{subsec:1.2};
in the second case, the pseudo-active effect is at work.

Generally, in a single-DoF system, semi-active vibration reduction is not possible because energy transfer to another degree of freedom cannot occur;
hence, only the pseudo-active effect is present.

\subsection{Serial two-DoF system}
\label{subsec:2.4}
\begin{figure} [hbt]
	\centering
	\includegraphics{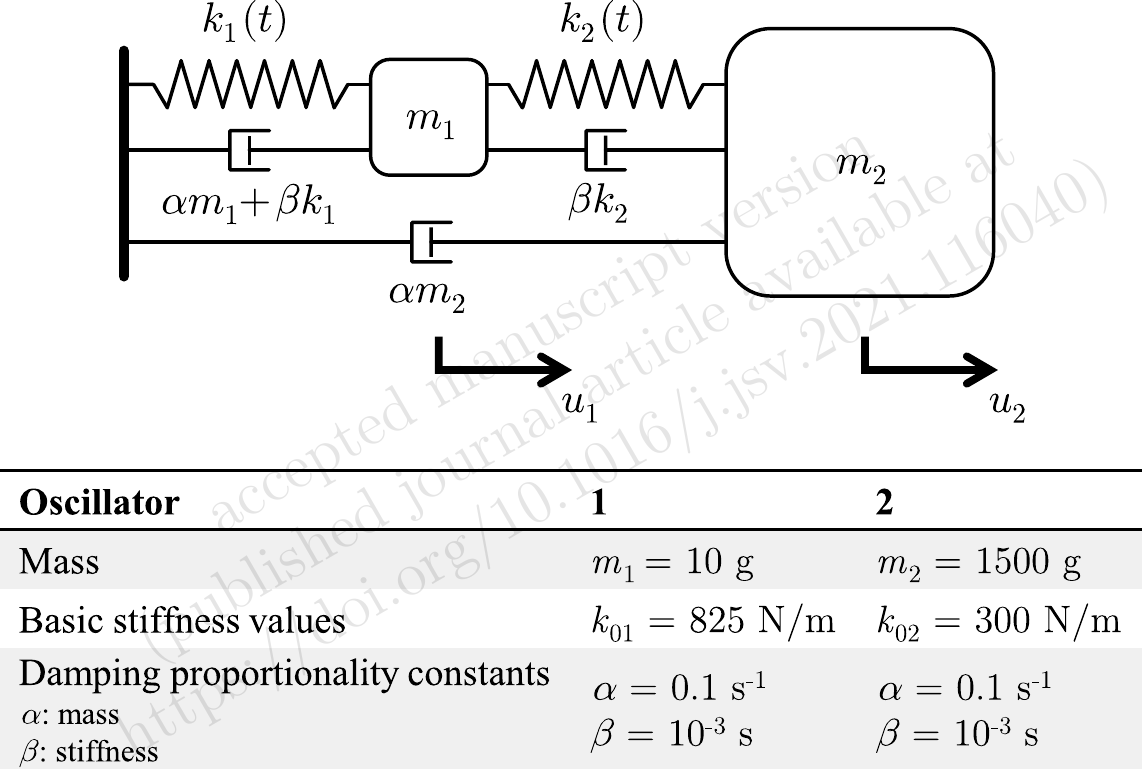}
	\caption{Serial two-DoF-system with variable stiffness}
	\label{fig:04}
\end{figure}

Now the example of a serial two-DoF spring-mass oscillator with variable stiffness is considered.
The system (\autoref{fig:04}) is appropriate for modeling the device in \cite{Ramaratnam2006}, which, as explained above, realizes two different stiffness values by changing the effective length of a spring.

The first spring (with stiffness $ k_1 $) in \autoref{fig:04} represents the portion of the spring used in \cite{Ramaratnam2006}, which can be rendered inactive;
the second spring (with stiffness $ k_2 $) represents the free portion.
The mass $ m_1 $, which is significantly smaller than $ m_2 $, represents the mass of the physical spring.
Mass- and stiffness-proportional damping is applied here with the same coefficients for both subsystems.

To directly reproduce the experiment presented in \cite{Ramaratnam2006}, $ k_1 $ must change between a given finite value to infinity, while $ k_2 $ remains at a constant value.
An alternative option, which reproduces the stiffness variation effect on the motion of mass $ m_2 $ by scaling both stiffness values, is presented later.

We define a low-stiffness state where
\begin{equation}
	k_1 = k_{01}, \quad
	k_2 = k_{02}
	\label{eq:06}
\end{equation}

In the high-stiffness state, the spring constants are scaled as follows:
\begin{equation}
	k_1 = \gamma_1 k_{01}, \quad
	k_2 = \gamma_2 k_{02}
	\label{eq:07}
\end{equation}

As will be shown later, owing to the relatively low value of mass $ m_1 $, the dynamics of the mass $ m_2 $ is essentially affected by the stiffness $ k $ of the serial arrangement of the two springs
\begin{equation}
	k = \dfrac{k_1 k_2}{k_1 + k_2}
	\label{eq:08}
\end{equation}

The low- and high-stiffness states are, respectively,
\begin{equation}
	k_l = \dfrac{k_{01} k_{02}} {k_{01} + k_{02}}, \quad
	k_h = \dfrac{\gamma_1 \gamma_2 k_{01} k_{02}} {\gamma_1 k_{01} + \gamma_2 k_{02}}
	\label{eq:09}
\end{equation}

The aforementioned option to reproduce the experiment in \cite{Ramaratnam2006} (hereafter referred to as the \textit{local} case) requires
\begin{equation}
	\gamma_1 \rightarrow \infty, \quad
	\gamma_2 = 1
	\label{eq:10}
\end{equation}
which, with \autoref{eq:09} leads to
\begin{equation}
	k_l = \dfrac{k_{01} k_{02}} {k_{01} + k_{02}}, \quad
	k_h = k_{02}
	\label{eq:11}
\end{equation}

By equating with the values of $ k_l $ and $ k_h$ reported in \cite{Ramaratnam2006}
\begin{equation}
	k_l = 220\,\mathrm{N/m}, \quad
	k_h = 300\,\mathrm{N/m}
	\label{eq:12}
\end{equation}
the basic stiffness values $ k_{01} $ and $ k_{02} $ specified in \autoref{fig:04} are obtained.

As mentioned, the effect on the motion of mass $ m_2 $ can also be reproduced using the same scaling factor for both springs (\textit{global} case).
With
\begin{equation}
	\gamma_1 = \gamma_2 = \gamma
	\label{eq:13}
\end{equation}
equation~\autoref{eq:09} writes
\begin{equation}
	k_l = \dfrac{k_{01} k_{02}} {k_{01} + k_{02}}, \quad
	k_h = \gamma \dfrac{k_{01} k_{02}} {k_{01} + k_{02}}
	\label{eq:14}
\end{equation}

By using the already determined basic stiffness values and \autoref{eq:14}, the value $ \gamma = 1.3636 $ is obtained.

\begin{figure} [hbt]
	\centering
	\captionsetup{justification=centering}
	\includegraphics{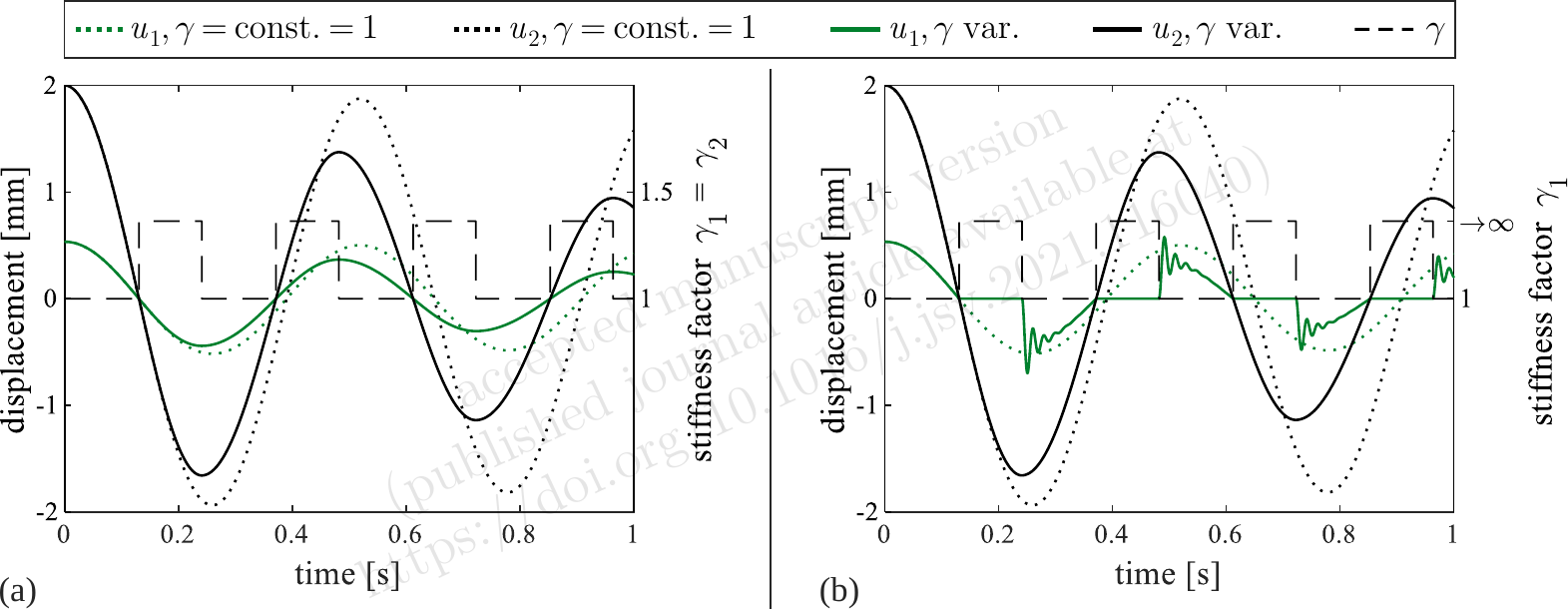}
	\caption{Representation of the experimental setup in \cite{Ramaratnam2006} with aid of a serial two-DoF system;\\(a) global case; (b) local case}
	\label{fig:05}
\end{figure}

The initial conditions are given by imposing a static displacement of the second degree of freedom $u_2$ (no velocity) and leaving the first degree of freedom $u_1$ unloaded.
The control logic is defined by \autoref{eq:02} with the observation variable $c=u_2$, and the stiffness scaling factors in the high-stiffness phase given by \autoref{eq:10} and \autoref{eq:13}, respectively.
The resulting displacement curves of the two variants are shown in \autoref{fig:05}, each compared with the non-modulated structure (dotted lines).

In the global case (a), both degrees of freedom resemble -- like the green curve in \autoref{fig:03} -- a damped sinusoidal curve.

The real law of the curves is not a damped sine but a piecewise sinusoidal curve (\autoref{fig:03}) or a piecewise damped sinusoidal curve (\autoref{fig:05}a).
We now consider the green curve in \autoref{fig:03}.
Because the system is undamped and the stiffness is piecewise constant, the curve will be piecewise harmonic.
In the theoretical case in which the switching stiffness control logic is applied exactly at the extremal and zero-crossing points, the low- and high-stiffness phases will respectively have the following lengths:
\begin{equation}
	P_l = \frac{\pi}{2} \cdot \sqrt{\frac{m_1}{k_1}}, \quad
	P_h = \frac{\pi}{2} \cdot \sqrt{\frac{m_1}{k_1+k_2}}
	\label{eq:15}
\end{equation}

Any of the phases is denoted as a \textit{quarter cycle} in the following.
In the first quarter cycle in \autoref{fig:03}, there is no difference with respect to the free oscillation of an undamped system with a finite displacement and zero velocity as the initial condition.
At the first zero-crossing, the stiffness is increased without affecting the vibration energy (as shown later).
In the second quarter cycle, the spring, loaded with the same energy present in the system at the initial time, reaches a lower amplitude owing to the higher stiffness.
In addition, we observe a shortening of the quarter cycle owing to the higher natural frequency.
At the end of the second quarter cycle, the stiffness is decreased, and energy is subtracted from the system.
The quarter cycle resumes its original length, and the velocity at the next zero-crossing decreases (owing to the loss of vibration energy), which can be seen in the decrease in the curve slope.
The fourth quarter cycle behaves essentially like the second cycle.

In a similar way, it could be shown that the curves in \autoref{fig:05}a are -- due to the presence of damping -- piecewise damped sinusoidal curves.

In the local case (\autoref{fig:05}b), the first degree of freedom remains blocked starting from a zero-crossing point, whereas the absolute displacement value of the second degree of freedom increases.
At each extremal point of the $ u_2 $ curve, the first spring is released, which leads to clearly visible high-frequency oscillations of the first degree of freedom.
Although the course of the displacement of the first degree of freedom differs significantly between the global and local case, the amplitude reduction for the second degree of freedom shows no appreciable differences between the two cases.

The undamped equations of motion, valid for both variants, are as follows
\begin{align}
\begin{split}
	m_1 \ddot{u}_1 + (k_1 + k_2) u_1 - k_2 u_2 &= 0\\
	m_2 \ddot{u}_2 - k_2 u_1 + k_2 u_2 &= 0
	\label{eq:16}
\end{split}
\end{align}

Owing to the low value of $ m_1 $, the inertial term can be neglected in the first equation, unless the accelerations are particularly high.
The following relationship is then obtained:
\begin{equation}
	u_1 = \dfrac{k_2}{k_1 + k_2} u_2
	\label{eq:17}
\end{equation}

The two displacements are proportional to each other.
The mass $ m_1 $ just ``follows'' $ m_2 $ in a quasi-static equilibrium state.
The system can be reduced to a single DoF, with the equation (obtained from the second equation of \autoref{eq:16} and \autoref{eq:17})
\begin{equation}
	m_2 \ddot{u}_2 + \dfrac{k_1 k_2}{k_1 + k_2} u_2 = 0
	\label{eq:18}
\end{equation}

Because the stiffness term of this equation yields, for both the local and global case, the same value in the high- and low-stiffness phase (see \autoref{eq:11} and \autoref{eq:14}), the curve for $ u_2 $ is approximately the same in \autoref{fig:05}a and b (due to the assumption $m_1 \ll m_2$).

In the global case, the curves differ only by a scaling factor, which is in accordance with \autoref{eq:17} and confirms the hypothesis that $ m_1 $ is negligible.
Because the spring stiffness values are scaled by the same factor, the proportionality coefficient in \autoref{eq:17} does not change.
According to \autoref{eq:14}, at the beginning of the observation, the system equation~\autoref{eq:18} can be written as
\begin{equation}
	m_2 \ddot{u}_2 + k_l u_2 = 0
	\label{eq:19}
\end{equation}

After one quarter cycle, the stiffness increases and the equation becomes
\begin{equation}
	m_2 \ddot{u}_2 + k_h u_2 = 0
	\label{eq:20}
\end{equation}
with the values for $ k_l $ and $ k_h $ given in \autoref{eq:14}.

The local case does not deviate from the global case in the first quarter cycle.
At the zero crossing, $ k_1 $ is set to an infinite value by imposing $ u_1 = 0 $.
The system is again reduced to a single-DoF system, with the equation (obtained from \autoref{eq:16}, second equation, for $ u_1 = 0 $)
\begin{equation}
	m_2 \ddot{u}_2 + k_2 u_2 = 0
	\label{eq:21}
\end{equation}
which, according to \autoref{eq:11} with $ k_2 = \mathrm{const.} = k_{02} $, provides equation~\autoref{eq:20}.

At the extremal point, spring~1 is released, and its stiffness is instantaneously set back to $ k_{01} $.
The force acting on $ m_1 $ suddenly increases from zero to a finite value, and the resulting acceleration renders the inertial term in \autoref{eq:16} no longer negligible.
This is clearly shown by the $ u_1 $ curve, which no longer follows the proportionality law \autoref{eq:17}.
The mass oscillates about the quasi-static equilibrium position given by the right-hand side of \autoref{eq:17}.
Toward the end of the quarter cycle, however, this effect swings out, and the curve again fulfills \autoref{eq:18}.

Because the inertial forces generated by $ m_1 $ cannot be neglected in this quarter cycle, equation~\autoref{eq:18} -- strictly speaking -- is no longer valid.
The $ u_2 $ curve, however, does not show any significant difference with respect to the global case.
Evidently, the fast oscillating inertial forces of $ m_1 $ do not affect the motion of the much larger mass $ m_2 $ to an appreciable extent.

The (total) vibration energy $ E $ is the sum of the elastic potential energy $ U $ and the kinetic energy $ T $.
The potential energy is given in both cases by the sum of the contributions of the two springs
\begin{equation}
	U = \frac{1}{2} k_1 u_1^2 + \frac{1}{2} k_2 (u_1 - u_2)^2
	\label{eq:22}
\end{equation}

As shown for the single-DoF system, the stepwise stiffness reductions, which occur at finite displacements, imply stepwise changes in the potential energy.
These energy changes are achieved by the stiffness variation device (extracted energy).
For the stiffness increase points, which occur at zero displacements and therefore do not alter the potential energy, no external work is required to balance changes in potential energy.

From \autoref{eq:22}, the expression for the extracted energy at the $i$th stiffness reduction point $ e_i $ can be written as
\begin{equation}
	\Delta E_{ai} = \frac{1}{2} \Delta k_1 u_1^2(e_i) + 
					\frac{1}{2} \Delta k_2(u_1(e_i) - u_2(e_i))^2
	\label{eq:23}
\end{equation}
where the negative factors $ \Delta k_j $ express the drop of the spring constants
\begin{equation}
	\Delta k_j = k_{0j} (1-\gamma_j)
	\label{eq:24}
\end{equation}

The vibration energy loss over the half cycle between $ e_i $ and the next extremal point $ e_{i+1} $, computed as the difference of the left limits of potential energies
\begin{equation}
	\Delta E_i = U(e_{i+1}^-) - U(e_i^-)
	\label{eq:25}
\end{equation}
corresponds, since the system is conservative, to the extracted energy.
In the general case in which the system's internal damping is considered, the vibration energy loss is expressed as the sum of the extracted energy and the energy dissipated by damping
\begin{equation}
	\Delta E_i = \Delta E_{ai} + \Delta E_{pi}
	\label{eq:26}
\end{equation}

Dividing by the vibration energy of the system at time $ e_i $
\begin{equation}
	\dfrac{\Delta E_i}{E(e_i)} = 
	\dfrac{\Delta E_{ai}}{E(e_i)} + 
	\dfrac{\Delta E_{pi}}{E(e_i)}
	\label{eq:27}
\end{equation}
renders the energy balance independent of the choice of the half cycle because the damping losses are proportional to the available vibration energy.
The same holds true for the extracted energy because it depends, according to \autoref{eq:23}, on the square of the displacement amplitudes.
The numerator and denominator of the fractions in \autoref{eq:27} change with time, but the quotients remain constant.
To avoid overloading the notation, we omit the dependence on $ i $ and $ e_i $ in the symbols of the energy loss rates in the following equations.

As long as relationship \autoref{eq:17} holds true, the system's potential energy \autoref{eq:22} can be expressed as a function of $ u_2 $.
This generally applies to the global case as well as during the high-stiffness phase of the local case.
In the high-stiffness phase, it holds
\begin{equation}
	U_h = \dfrac{1}{2} k_h u_2^2
	\label{eq:28}
\end{equation}
and \autoref{eq:25} can then be written as
\begin{equation}
	\Delta E_i = \dfrac{1}{2} k_h (u_2^2 (e_{i+1}^-) - u_2^2 (e_i^-))
	\label{eq:29}
\end{equation}

Assuming the same value of the function $ u_2 $ for both cases, it can be concluded that the local stiffness modulation leads to the same energy loss rate $ \Delta E / E $ as the global stiffness modulation.

In the local case, however, the extracted energy is zero, according to \autoref{eq:23}, because $ u_1 $ and $ \Delta k_2 $ are zero.
The damping energy rate in \autoref{eq:27} is defined to be equal to the energy loss rate of the system without stiffness variations, which is equal for the two cases.
Hence, a new term must replace the extracted energy in the half-cycle balance:
\begin{equation}
	\dfrac{\Delta E}{E} = 
	\dfrac{\Delta E_s}{E} + 
	\dfrac{\Delta E_p}{E}
	\label{eq:30}
\end{equation}

The new term $ \Delta E_s/E $ represents the semi-active effect.
In both cases, the same share of energy is periodically extracted from the motion of the second degree of freedom
\begin{equation}
	\dfrac{\Delta E_s}{E} = 
	\dfrac{\Delta E_a}{E}
	\label{eq:31}
\end{equation}
with the difference that, in the local case, instead of being extracted by the actuator, this energy is transferred to the high-frequency motion of the first degree of freedom and dissipated, similarly to the case in \autoref{fig:02}.

In the global case, with \autoref{eq:28} and the following expression of potential energy for the low-stiffness phase
\begin{equation}
	U_l = \dfrac{1}{2} k_l u_2^2
	\label{eq:32}
\end{equation}
the extracted energy \autoref{eq:23} can be written as
\begin{equation}
	\Delta E_{ai} = U_l - U_h = \dfrac{1}{2} (k_l - k_h) u_2^2 (e_i)
	\label{eq:33}
\end{equation}

Finally, the energy rates from \autoref{eq:31} can be expressed as 
\begin{equation}
	\dfrac{\Delta E_{a}}{E} = \dfrac{\Delta E_s}{E} =
	\dfrac{k_l - k_h}{k_h} = \dfrac{k_l}{k_h} - 1
	\label{eq:34}
\end{equation}

Note that in this energy analysis, we neglected the small amount of external, negative work to bring the velocity of mass $ m_1 $ to zero at the stiffness increase points (kinks in the $ u_1 $ displacement curves).

The presented example helped us to define the different contributions to vibration reduction by separating the pseudo-active and semi-active components.
In the general case, all the considered contributions coexist:
\begin{equation}
	\dfrac{\Delta E}{E} = 
	\dfrac{\Delta E_a}{E} + \dfrac{\Delta E_s}{E} + \dfrac{\Delta E_p}{E}
	\label{eq:35}
\end{equation}

The lost energy as a function of time
\begin{equation}
	L(t) = E(t) - E(0)
	\label{eq:36}
\end{equation}
contains the three components resulting from the analyzed effects
\begin{equation}
	L(t) = L_a(t) + L_s(t) + L_p(t)
	\label{eq:37}
\end{equation}

The extracted energy $ L_a(t) $ is a piecewise constant function given by
\begin{equation}
	L_a(t) = \sum_i \Delta E_{ai} H(e_i)
	\label{eq:38}
\end{equation}
with $ H $ as the Heaviside function and $\Delta E_{ai}$ resulting from \autoref{eq:23}.

The difference between the lost energy and extracted energy is due to internal dissipation.
It consists of the above introduced terms $ L_s(t) $ and $ L_p(t) $, which are continuous functions of time, and can be analytically expressed as
\begin{equation}
\delimitershortfall -0.1pt
	L_s(t) + L_p(t) = 
	 -\int_0^t
		\dot{u}_1^2(\tau) \left(\alpha m_1 + \beta (k_1+k_2)\right) + 
		\dot{u}_2^2(\tau) (\alpha m_2 + \beta k_2) - 
		2 \dot{u}_1(\tau) \dot{u}_2(\tau) \beta k_2 \,d\tau
	\label{eq:39}
\end{equation}
with varying $ k_1 $ and $ k_2 $ as a function of $ \gamma $ \autoref{eq:02}.

\begin{figure} [hbt]
	\centering
	\captionsetup{justification=centering}
	\includegraphics{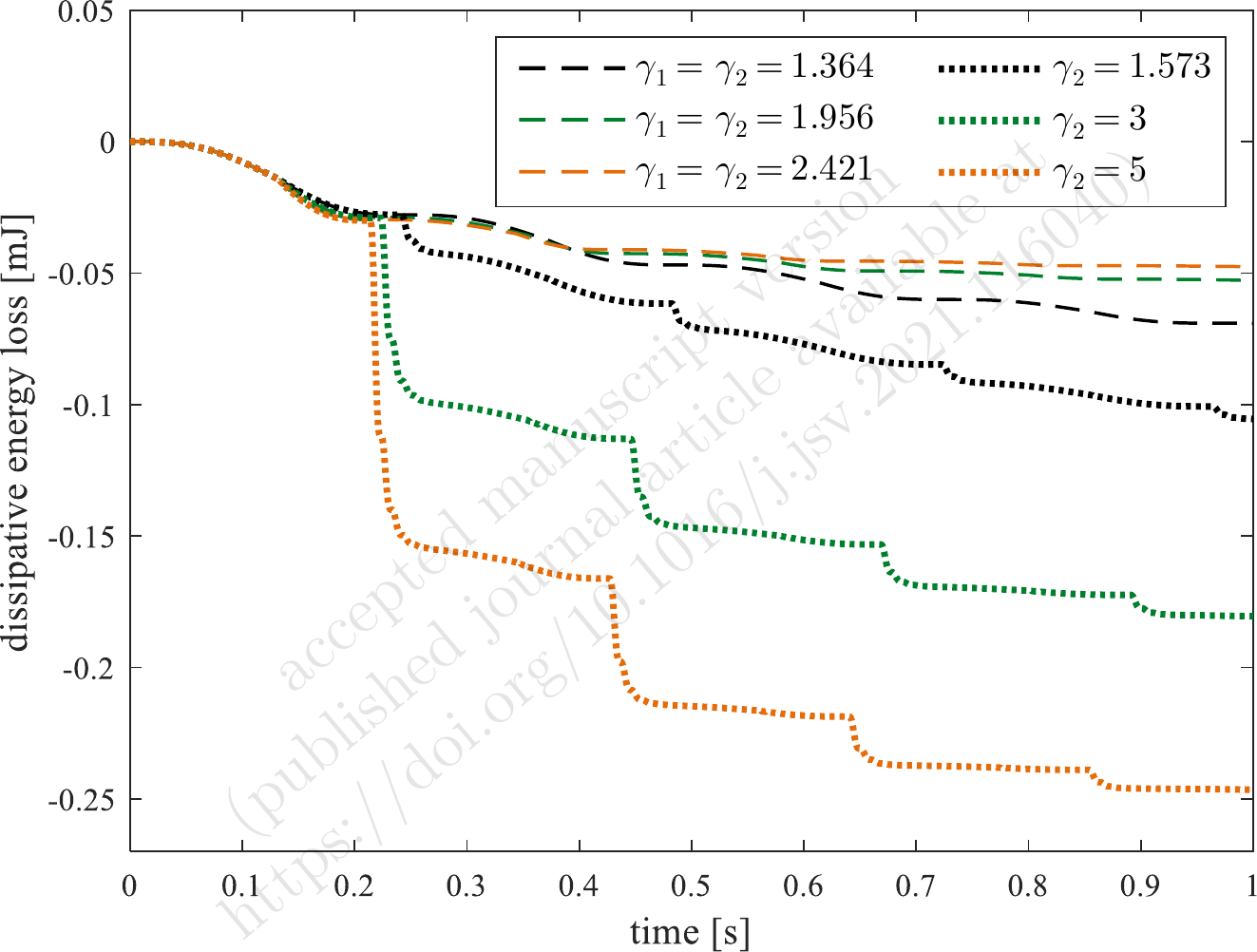}
	\caption{Dissipative energy loss $ L_s(t) + L_p(t) $;\\
			 dashed lines: global stiffness scaling;
			 dotted lines: local stiffness scaling of spring 2}
	\label{fig:06}
\end{figure}

The diagram in \autoref{fig:06} shows the dissipative energy loss $ L_s(t) + L_p(t) $ for different local and global stiffness variations.
Other than in the previous example, spring~2 is affected by the stiffness variation in the local case, which implies that $ L_a $ is now also present.
In the global case (dashed lines) the curves only represent the contribution of the passive term $ L_p(t) $, while the curves of the local case (dotted lines) additionally include the semi-active contribution $ L_s(t) $.
This is shown later in more detail.
The stiffness scaling factors of the three local examples are chosen in such a way that the total dissipated energy of the system per cycle is pairwise identical to that of the global examples.
Curves of the same color refer to example pairs with the same energy loss.
Because the passive term is assumed to be equal for both cases, the difference between the two curves of the same color shows the semi-active component $ L_s(t) $.

With higher values of the factor $ \gamma $, the dissipative energy loss clearly increases in the local case, mainly because of the semi-active effect.
In the global case, the passive component shows a reverse trend with respect to the stiffness scaling factor.
This can be attributed to the fact that the percentage reduction in vibration energy by damping is proportional to the vibration energy itself.
If another vibration reduction mechanism is at work (pseudo-active in the present case), it reduces the vibration amplitude and therefore indirectly diminishes the \textit{absolute} energy dissipated by inherent damping.

\subsection{Modal analysis}
\label{subsec:2.5}
In the analysis of the system with coupleable oscillators (\autoref{fig:02}), the energy transfer between the two oscillators is evident.
In a serial system, such as the one examined in the previous section, the two oscillators are always coupled;
hence, the vibration energy cannot be directly assigned to the individual physical degrees of freedom.
However, we were able to study the stiffness modulation effect energetically because of the relatively small mass of one degree of freedom, but this cannot be expected to represent the general case.
In the modal space, decoupled oscillators are present (which recalls the first example), and energy can be split over the modal degrees of freedom.
In this section, we therefore analyze the stiffness modulation effect by modal analysis.
The analysis is specifically performed on the serial two-DoF system, but can easily be extended to the general case of a discrete system with a larger number of degrees of freedom or a continuous structure.

The stiffness matrix $ \mathbf{K} $ and mass matrix $ \mathbf{M} $ of the system can be calculated as follows:
\begin{equation}
	\mathbf{K} = \begin{bmatrix}
		k_{01} + k_{02} &	-k_{02}\\
		-k_{02}			&	 k_{02}
	\end{bmatrix}, \quad
	\mathbf{M} = \begin{bmatrix}
		m_1 	& 	0\\
		0		&	m_2
	\end{bmatrix}
	\label{eq:40}
\end{equation}

Damping is assumed to be proportional to mass and stiffness, and the damping matrix $ \mathbf{C} $ thus results from the linear combination of $ \mathbf{K} $ and $ \mathbf{M} $
\begin{equation}
	\mathbf{C} = \alpha \mathbf{M} + \beta \mathbf{K}
	\label{eq:41}
\end{equation}

The solution of the eigenvalue problem
\begin{equation}
	\mathbf{K} \boldsymbol{\unslant\varphi} = \lambda \mathbf{M} \boldsymbol{\unslant\varphi}
	\label{eq:42}
\end{equation}
provides the system's eigenvalues $ \lambda_1, \lambda_2 $ and the eigenvectors $ \boldsymbol{\unslant\varphi}_1, \boldsymbol{\unslant\varphi}_2 $.
It holds
\begin{equation}
	\lambda_i = \omega_i^2,\quad
	i = 1,2
	\label{eq:43}
\end{equation}
with $ \omega_i $ as the system's eigenfrequencies.

The transformation into the modal space follows the relationship
\begin{equation}
	\mathbf{u = \Phi}^\mathrm{T} \mathbf{q}
	\label{eq:44}
\end{equation}
with $ \mathbf{u} $ as the spatial displacements, $ \mathbf{q} $ as the modal displacements, and
\begin{equation}
	\mathbf{\Phi} = [\,\boldsymbol{\unslant\varphi}_1 \;|\; \boldsymbol{\unslant\varphi}_2\,]
	\label{eq:45}
\end{equation}
as the modal matrix.
The modal displacement fulfills the modal system equation
\begin{equation}
	\tilde{\mathbf{M}} \ddot{\mathbf{q}} + 
	\tilde{\mathbf{C}} \dot{\mathbf{q}} + 
	\tilde{\mathbf{K}} \mathbf{q} = 
	\tilde{\mathbf{f}}
	\label{eq:46}
\end{equation}
with the modal system matrices:
\begin{align}
\begin{split}
	\tilde{\mathbf{M}} &= \mathbf{\Phi}^\mathrm{T} \mathbf{M} \mathbf{\Phi}\\
	\tilde{\mathbf{C}} &= \mathbf{\Phi}^\mathrm{T} \mathbf{C} \mathbf{\Phi}\\
	\tilde{\mathbf{K}} &= \mathbf{\Phi}^\mathrm{T} \mathbf{K} \mathbf{\Phi}
	\label{eq:47}
\end{split}
\end{align}

The diagonal entries $ \tilde{m}_i $, $ \tilde{c}_i $, and $ \tilde{k}_i $ of the modal matrices supply the system parameters of the decoupled oscillators.
Owing to the mass and stiffness proportionality, not only $ \tilde{\mathbf{M}} $ and $ \tilde{\mathbf{K}} $ but also $ \tilde{\mathbf{C}} $ is a diagonal matrix.

First, we consider how a global stiffness change affects the modal quantities.
In this case, all the elements of the stiffness matrix $ \mathbf{K} $ are scaled by the same factor $ \gamma_1 = \gamma_2 = \gamma $:
\begin{equation}
	\mathbf{K}_h = \gamma \mathbf{K}
	\label{eq:48}
\end{equation}

The eigenvalue problem of the new system
\begin{equation}
	\mathbf{K}_h \boldsymbol{\unslant\varphi} = 
	\bar{\lambda} \mathbf{M} \boldsymbol{\unslant\varphi}
	\Rightarrow
	\gamma \mathbf{K} \boldsymbol{\unslant\varphi} = 
	\bar{\lambda} \mathbf{M} \boldsymbol{\unslant\varphi}
	\label{eq:49}
\end{equation}
is fulfilled by the same eigenvectors of \autoref{eq:42} and by the eigenvalues
\begin{equation}
	\bar{\lambda}_i = \gamma \lambda_i
	\label{eq:50}
\end{equation}

The modal transformation does not change.
The modal spring-mass systems experience the same proportional change in stiffness as the physical one, while the modal masses remain unchanged.
Owing to the change in the modal stiffness, the energy also changes in each modal oscillator separately (by the same factor $ \gamma $), but no energy transfer occurs between the oscillators.
As previously observed in the single-DoF case, a global stiffness variation can therefore only lead to a pseudo-active and not a semi-active vibration suppression.

Now, we consider the case of a local change in stiffness, realized by scaling the second spring stiffness by the factor $ \gamma_2 $ and leaving $ k_1 $ unchanged.
In this case, the stiffness changes as follows:
\begin{equation}
	\mathbf{K}_h = \begin{bmatrix}
		k_{01}+ \gamma_2 \cdot k_{02} 	&	- \gamma_2 \cdot k_{02}\\
		- \gamma_2 \cdot k_{02}		&	\gamma_2 \cdot k_{02}
	\end{bmatrix}
	\label{eq:51}
\end{equation}

As a result of the nonproportional modification of individual entries in the stiffness matrix, the modal basis $ \mathbf{\Phi} $ changes.
Consequently, at the moment of stiffness change, the vibration energy is redistributed among the modal oscillators, which in turn change their modal parameters.
Thus, the basis for vibration reduction through internal energy transfer is given.

\begin{figure} [p]
	\centering
	\includegraphics{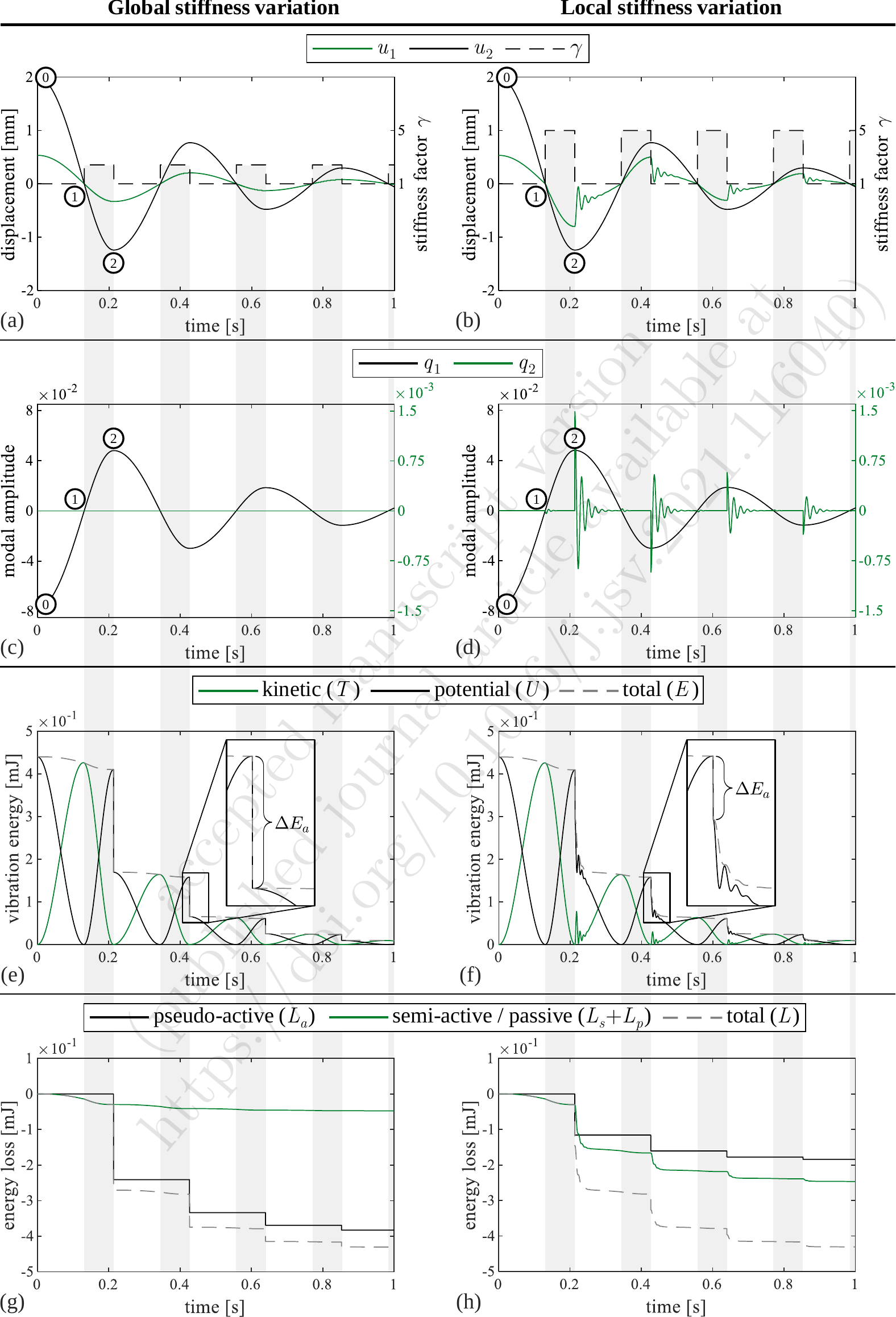}
	\caption{Results for global stiffness variation (left column) and local stiffness variation (right column);\\
	1\textsuperscript{st} row: physical displacement,
	2\textsuperscript{nd} row: modal amplitudes,
	3\textsuperscript{rd} row: vibration energy,\\
	4\textsuperscript{th} row: energy loss, separated into pseudo-active and semi-active + passive component}
	\label{fig:07}
\end{figure}

The described difference will now be illustrated by means of an exemplary calculation with the data given in \autoref{tab:1}.
The stiffness scaling factors correspond to the orange curve pair in \autoref{fig:06}, and result in the same total energy loss over the observed time period.
Analogously to the previous examples, the stiffness is switched according to \autoref{eq:02} with the first modal amplitude $ q_1 $ as observation variable $ c $ (see \autoref{fig:07}c and d).
The previously used formulation with $c=u_2$ yields similar switching points (see \autoref{fig:07}a and b) but should be omitted here for the purpose of explanation.

\begin{table}[ht]
	\centering
	\footnotesize
	\begin{tabularx}{16cm}{p{3.2cm}| X| X| X} & 
		\centering \textbf{Unchanged} &
		\centering \textbf{Global change} &
		\centering \textbf{Local change} \tabularnewline
		\hline
		\rule{0pt}{3ex}%
		\textbf{Stiffness scaling factors} &
		\centering -- &
		\centering $ \gamma_1 = \gamma_2 = \gamma = 2.421 $ &
		\centering $ \gamma_2=5 $ \tabularnewline[5pt]

		\rule{0pt}{3ex}%
		\textbf{Modal mass matrix} &
		$\tilde{\mathbf{M}} = $ \scriptsize $\begin{bmatrix}
			1 & 0\\0 & 1\end{bmatrix}$ &
		$\tilde{\mathbf{M}} = $ \scriptsize $\begin{bmatrix}
			1 & 0\\0 & 1\end{bmatrix}$ &
		$\tilde{\mathbf{M}} = $ \scriptsize $\begin{bmatrix}
			1 & 0\\0 & 1\end{bmatrix}$\\

		\rule{0pt}{5ex}%
		\textbf{Modal damping matrix} &
		$\tilde{\mathbf{C}} = $ \scriptsize $\begin{bmatrix}
			0.247 & 0\\0 & 112.653\end{bmatrix}$ &
		$\tilde{\mathbf{C}} = $ \scriptsize $\begin{bmatrix}
			0.455 & 0\\0 & 272.592\end{bmatrix}$ &
		$\tilde{\mathbf{C}} = $ \scriptsize $\begin{bmatrix}
			0.454 & 0\\0 & 233.246\end{bmatrix}$ \\

		\rule{0pt}{5ex}%
		\textbf{Modal stiffness matrix}	&
		$\tilde{\mathbf{K}} = $ \scriptsize $\begin{bmatrix}
			146.597 & 0\\0 & 1.126\cdot10^5\end{bmatrix}$ &
		$\tilde{\mathbf{K}} = $ \scriptsize $\begin{bmatrix}
			354.912 & 0\\0 & 2.725\cdot10^5\end{bmatrix}$ &
		$\tilde{\mathbf{K}} = $ \scriptsize $\begin{bmatrix}
			353.855 & 0\\0 & 2.332\cdot10^5\end{bmatrix}$ \\

		\rule{0pt}{5ex}%
		\textbf{Modal matrix} &
		$\mathbf{\Phi} = $ \scriptsize $\begin{bmatrix}
			-6.893 & -316.153\\-25.814 & 0.563\end{bmatrix}$ &
		$\mathbf{\Phi} = $ \scriptsize $\begin{bmatrix}
			-6.893 & -316.153\\-25.814 & 0.563\end{bmatrix}$ &
		$\mathbf{\Phi} = $ \scriptsize $\begin{bmatrix}
			-16.66 & -315.789\\-25.784 & 1.36\end{bmatrix}$ \\
	\end{tabularx}
	\caption{Stiffness scaling factors and modal matrices}
	\label{tab:1}
\end{table}

At $ t = t_0 = 0 $ only the first mode is active and with zero velocity
\begin{equation}
	\mathbf{q}(t_0) = \begin{bmatrix} q_1(t_0) \\ 0	\end{bmatrix}, \quad
	\dot{\mathbf{q}}(t_0) = \begin{bmatrix} 0 \\ 0	\end{bmatrix}, \quad
	q_1(t_0) = -0.0775
	\label{eq:52}
\end{equation}

First, the \textit{global} case is considered.
At $ t = z_i $, both modal oscillators are at their static equilibrium position (see \autoref{fig:07}c, Point~1):
\begin{equation}
	\mathbf{q}(z_i) = \begin{bmatrix} 0 \\ 0	\end{bmatrix}, \quad
	\dot{\mathbf{q}}(z_i) = \begin{bmatrix} \dot{q}_1(z_i) \\ 0	\end{bmatrix}
	\label{eq:53}
\end{equation}

At this point, the stiffness matrix is scaled according to \autoref{eq:48}.
As mentioned before, this occurs without external work because the modal displacements are zero;
therefore, the potential energy does not change:
\begin{equation}
	U(z_i) = \dfrac{1}{2} \mathbf{q}^\mathrm{T} (z_i) \tilde{\mathbf{K}} \mathbf{q}(z_i) = 
	\dfrac{1}{2} \mathbf{q}^\mathrm{T} (z_i) \gamma \tilde{\mathbf{K}} \mathbf{q}(z_i) = 0
	\label{eq:54}
\end{equation}

In the following quarter cycle, the structure oscillates at a slightly higher frequency and lower amplitude than the initial configuration owing to the higher stiffness value.
This is only apparently a vibration reduction because the vibration energy in the system remains unchanged.
As already mentioned, no change occurs in the modal transformation and the second modal oscillator remains at rest (see \autoref{fig:07}c, green curve).
Thus, the energetic considerations made for the global case in \autoref{subsec:2.4} are strictly valid, and we now have (without the need for a simplifying hypothesis) a single-DoF system.
The energy loss in the semicycle starting at $ t = e_i^- $ (Point~2 in the graph) is given by \autoref{eq:26}, with the pseudo-active component expressed by (see \autoref{eq:33})
\begin{equation}
	\Delta E_a = \dfrac{1}{2}(1-\gamma) \mathbf{q}^\mathrm{T} (e_i) \tilde{\mathbf{K}} \mathbf{q}(e_i)
	\label{eq:55}
\end{equation}
and the passive term $ \Delta E_p $ owing to damping of the first mode.

In the \textit{local} case, as seen in the spatial analysis, the semi-active effect joins the pseudo-active effect.
At the stiffness variation points, the modal basis changes between $ \mathbf{\Phi}_l $ (low stiffness) and $ \mathbf{\Phi}_h $ (high stiffness) such that the modal displacements of the low-stiffness basis $ \mathbf{q}_l $ must be converted to the new space with modal displacements $ \mathbf{q}_h $, and vice versa:
\begin{equation}
	\mathbf{q}_h = \mathbf{\Phi}_h^{-1} \mathbf{\Phi}_l \mathbf{q}_l, \quad
	\mathbf{q}_l = \mathbf{\Phi}_l^{-1} \mathbf{\Phi}_h \mathbf{q}_h
	\label{eq:56}
\end{equation}

This change in the modal basis has two consequences: it changes the modal matrices (see \autoref{tab:1}), and therefore the parameters of the modal oscillators; and it rearranges the contributions of the single modes to the motion of the system.
The latter effect leads to the excitation of the second (modified) modal oscillator, as shown in \autoref{fig:07}d.
At the points where the stiffness is changed to its lower value, the second modal oscillator begins to oscillate at a high frequency, which recalls \autoref{fig:02} (the corresponding curves in row~2 are plotted on a separate axis with different scaling for better resolution).

As shown by previous examples and the energy analysis in \autoref{fig:07}, the energy transferred to the higher mode can dissipate much faster than in the first mode, which leads to a significant semi-active contribution to vibration reduction.
Before discussing the energy curves, it should be recalled that the stiffness factor $ \gamma $ is adjusted such that the total energy per cycle is identical for both methods.

The plots in \autoref{fig:07}e and f (3\textsuperscript{rd} row) show the potential, kinetic, and total energy as a function of time for both cases.
In the global case, between the extremal points, the total energy loss owing to system damping is minor.
The main part of the energy loss is provided by the steps at the extremal points owing to the extracted energy $ L_a $ reported in \autoref{fig:07}g, together with the passive component $ L_p $.
As mentioned, there is no excitation of the second modal oscillator (see \autoref{fig:07}c); therefore, $ L_s $ is equal to 0.

While examining the local case (\autoref{fig:07}, right-hand side), a clear difference with respect to the global case appears in the behavior of the curve between the extremal points.
The step in the energy loss is much smaller and is followed by a steep (but continuous) energy loss that extends over the time interval in which the high-frequency oscillations of the second mode occur.
The extracted energy (\autoref{fig:07}h) now represents a much smaller portion of the total energy loss.
The missing part of the loss is replaced by the additional energy loss occurring between the extremal points, which corresponds to $ \Delta E_s $.
As shown in \autoref{fig:07}h, the total internal energy dissipation $ L_s + L_p $  exceeds, for the case under consideration, its pseudo-active counterpart $ L_a $.

\begin{figure} [hbt]
	\centering
	\includegraphics{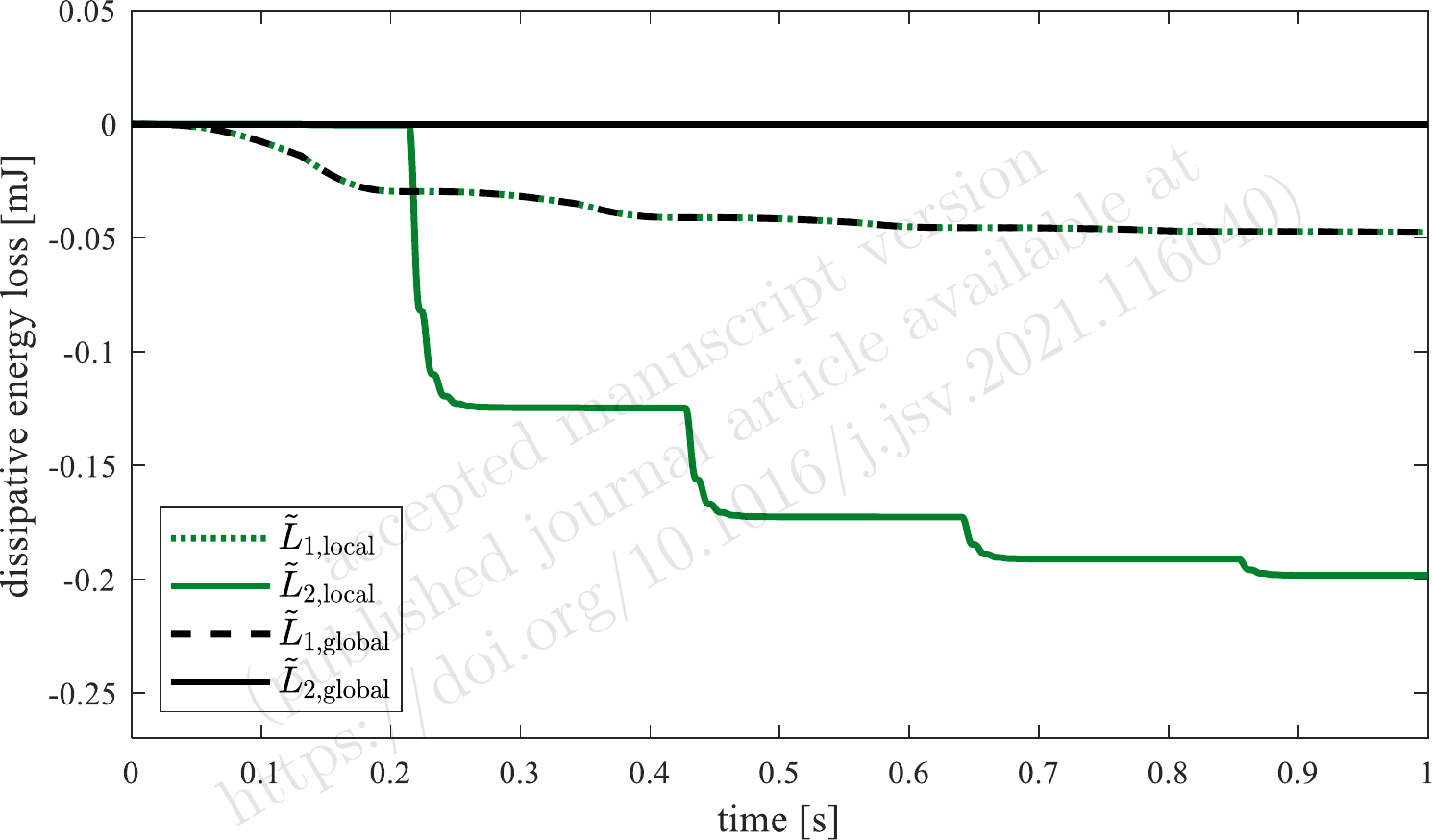}
	\caption{Dissipative energy loss, separated into modal components}
	\label{fig:08}
\end{figure}

In \autoref{fig:08}, the dissipative energy loss is shown separately for both modes and for the local and global case.
The curves represent the damped energy $ \tilde{L}_i $ related to mode $ i $:
\begin{equation}
	\tilde{L}_i (t,\gamma) = 
		-\int_0^t \tilde{c}_i (\gamma) \dot{q}_i^2(\tau) \,d\tau = 
		-\int_0^t \left(\alpha \tilde{m}_i + \beta \tilde{k}_i (\gamma) \right) \dot{q}_i^2(\tau) \,d\tau
	\label{eq:57}
\end{equation}

The dissipation related to mode~1 is equal for the two cases (dashed black and dotted green curve) and represents the passive component:
\begin{equation}
	\tilde{L}_{1,\mathrm{local}} = \tilde{L}_{1,\mathrm{global}} = L_p
	\label{eq:58}
\end{equation}

As expected, no dissipation occurs in the second mode in the global case (solid black curve) because the mode is not excited at all.
The dissipation associated with the second mode in the local case (solid green curve) constitutes the semi-active part $ L_s $:
\begin{equation}
	\tilde{L}_{2,\mathrm{local}} = L_s, \quad \tilde{L}_{2,\mathrm{global}} = 0
	\label{eq:59}
\end{equation}

\subsection{Forced oscillation}
\label{subsec:2.6}
Now, the case of a harmonic external force is considered.
The proof of the feasibility of the presented method is performed on the serial two-DoF system shown in \autoref{fig:04} with the stiffness amplification factors $ \gamma $ used in the previous section (see \autoref{tab:1}).

The analysis of the system without stiffness modification yields the natural frequencies 1.927\,Hz and 53.395\,Hz.
To determine the response behavior around the first natural frequency, an input force $ F(t) $ with an amplitude of 1\,N and a linearly increasing frequency between 1 and 3\,Hz is applied to mass $ m_2 $:
\begin{equation}
\begin{gathered}
	\delimitershortfall -0.1pt
	F(t) = F_{\mathrm{max}} \cdot \sin
			\left(2 \pi t \left(f_0 + t \dfrac{f_1-f_0}{2 t_1}\right) \right), \\
	t_1 = 100\,\mathrm{s}, \quad
	f_0 = 1\,\mathrm{Hz}, \quad
	f_1 = 3\,\mathrm{Hz}, \quad
	F_{\mathrm{max}} = 1\,\mathrm{N}
	\label{eq:60}
\end{gathered}
\end{equation}

\begin{figure} [hbt]
	\centering
	\includegraphics{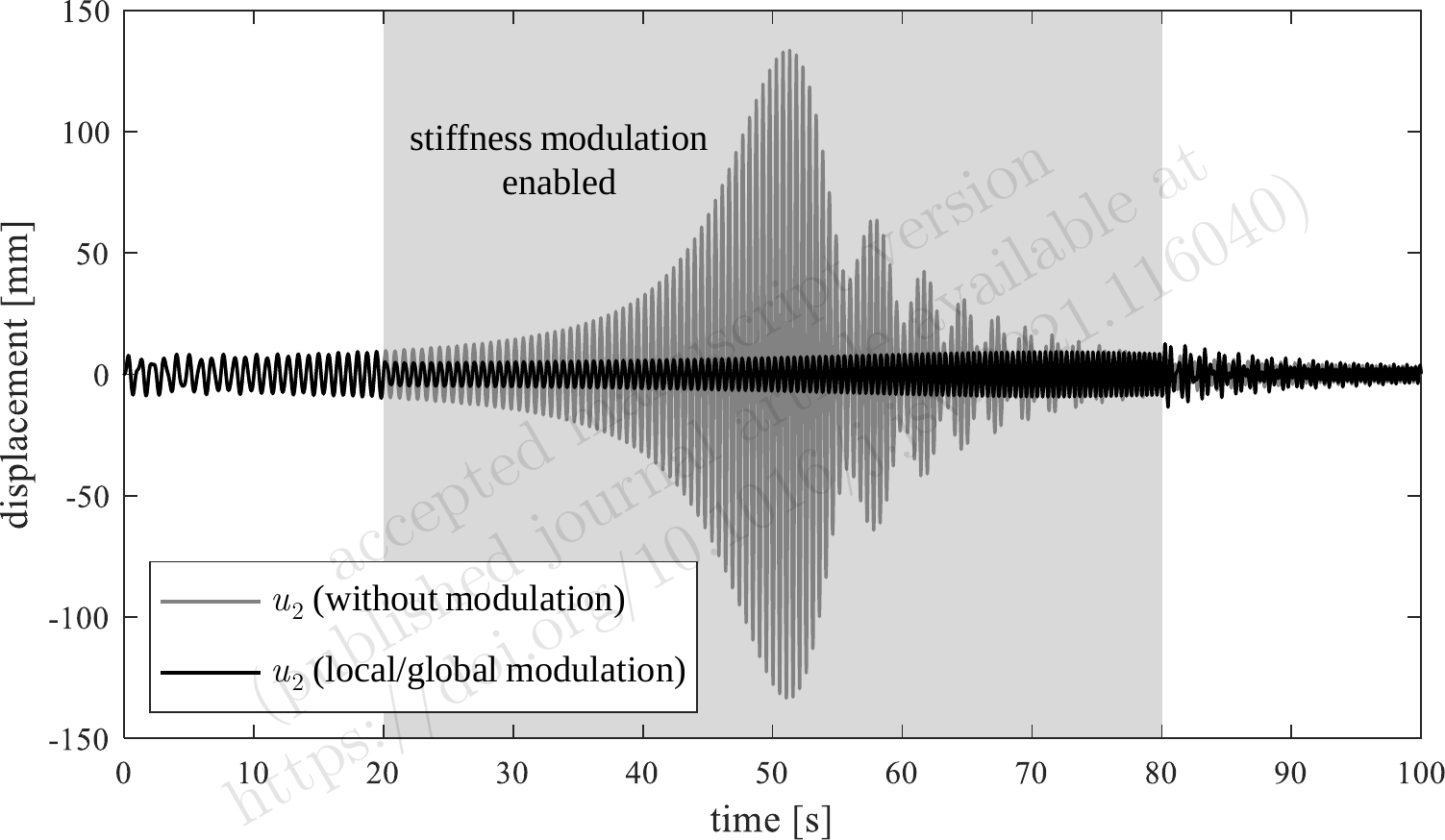}
	\caption{Displacement of $ u_2 $ due to a linear sweeping input force}
	\label{fig:09}
\end{figure}

The resulting displacement time graph of DoF $ u_2 $ is shown in \autoref{fig:09}.

The time response of $ u_2 $ without a change in stiffness (\autoref{fig:09}, gray curve) shows the typical resonance behavior in the neighborhood of the first eigenfrequency.
By modulating the stiffness in the same way as shown in \autoref{fig:07}, 1\textsuperscript{st} row (switch at zero-crossings / extremal points) this resonance amplification can be completely avoided (\autoref{fig:09}, black curve).
The stiffness modulation here operates only in the region indicated by the gray box.
Note that the black curve shows the time response of $ u_2 $ for both global and local stiffness variations.
As shown for the free oscillation, there are no deviations between the two cases for this degree of freedom with the selected parameter $ \gamma $ (for comparison see \autoref{fig:05} and \autoref{fig:07}).

In the free-vibration study it was shown that, while the total lost energy per cycle (due to the adjustment of the amplification factors) is the same in both cases, large differences between the local and global case exist in terms of the dissipated energy (by damping).
In the local case, the amount of dissipated energy is substantially higher than that in the global case.
The energy loss due to direct extraction at the stiffness variation interface is consequently higher in the global case.
The difference in dissipated energy between the two cases can be entirely attributed to the semi-active effect, since this is absent in the global case.

\begin{figure} [hbt]
	\centering
	\includegraphics[scale=0.92]{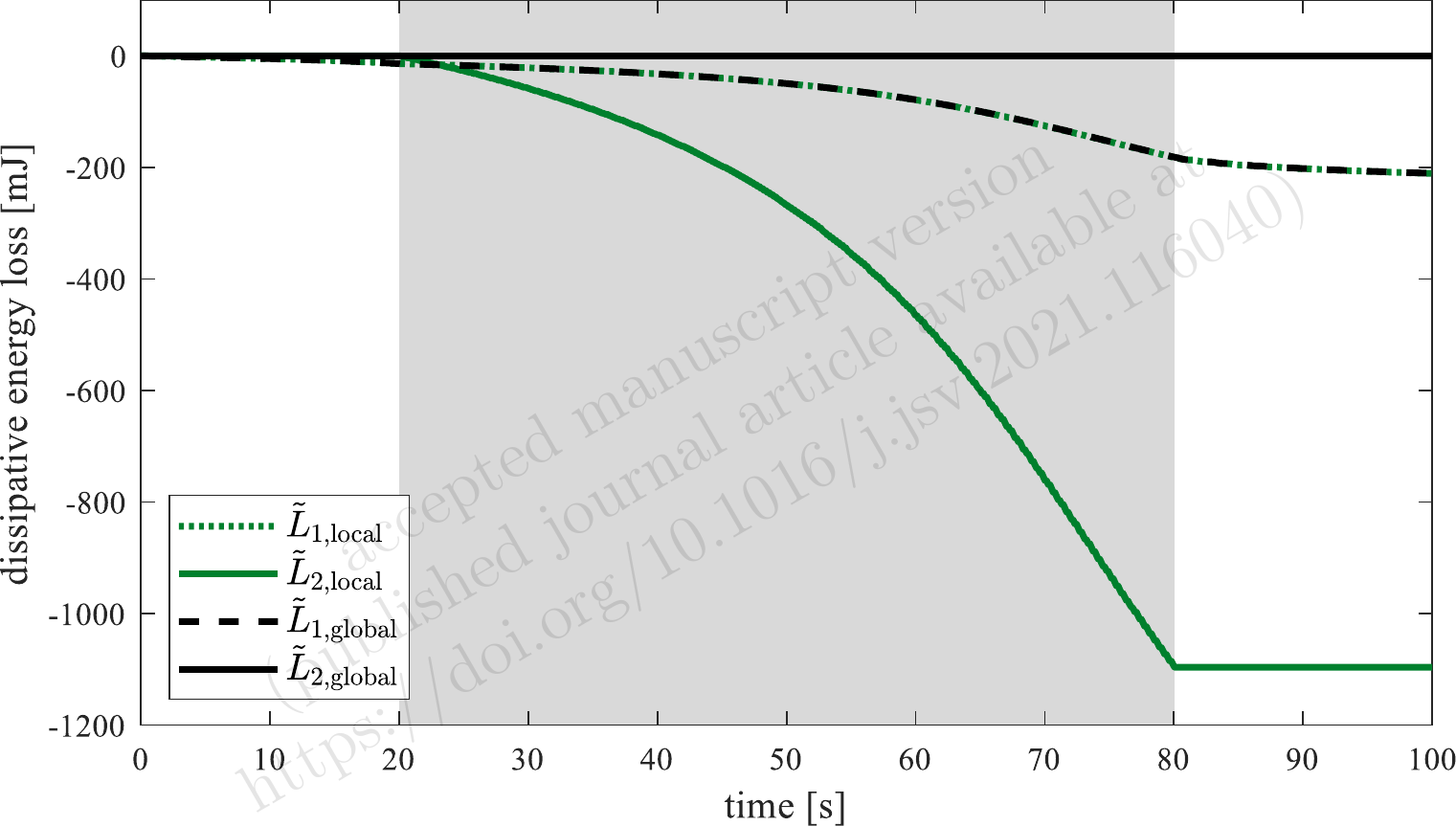}
	\caption{Dissipative energy loss, separated into modal components (sweep)}
	\label{fig:10}
\end{figure}

To extend this kind of analysis to the forced oscillation case, the dissipative energy loss in the two cases is analyzed for the individual modes according to \autoref{eq:57} (see \autoref{fig:10}).
The curve $ \tilde{L}_1 $ relative to the first mode is -- as previously stated in \autoref{eq:58} -- identical in both cases (dashed black and dotted green curve).
As far as the second mode is concerned, the applied force does not activate it because the excitation frequency lies far away from the second resonance.
In the global case, the dissipated energy of the second mode $ \tilde{L}_2 $ (solid black curve) remains negligible (as in the free vibration case; see \autoref{fig:08}); it is evident that no excitation of the second mode by energy transfer occurs.
Thus, it can be concluded that there is no semi-active effect in the global case (see also \autoref{eq:59}).
In the case of local stiffness modulation, a large amount of energy is dissipated by the second mode (solid green curve) during the stiffness variation phase, which shows the semi-active component.

\begin{figure} [hbt]
	\centering
	\includegraphics[scale=0.92]{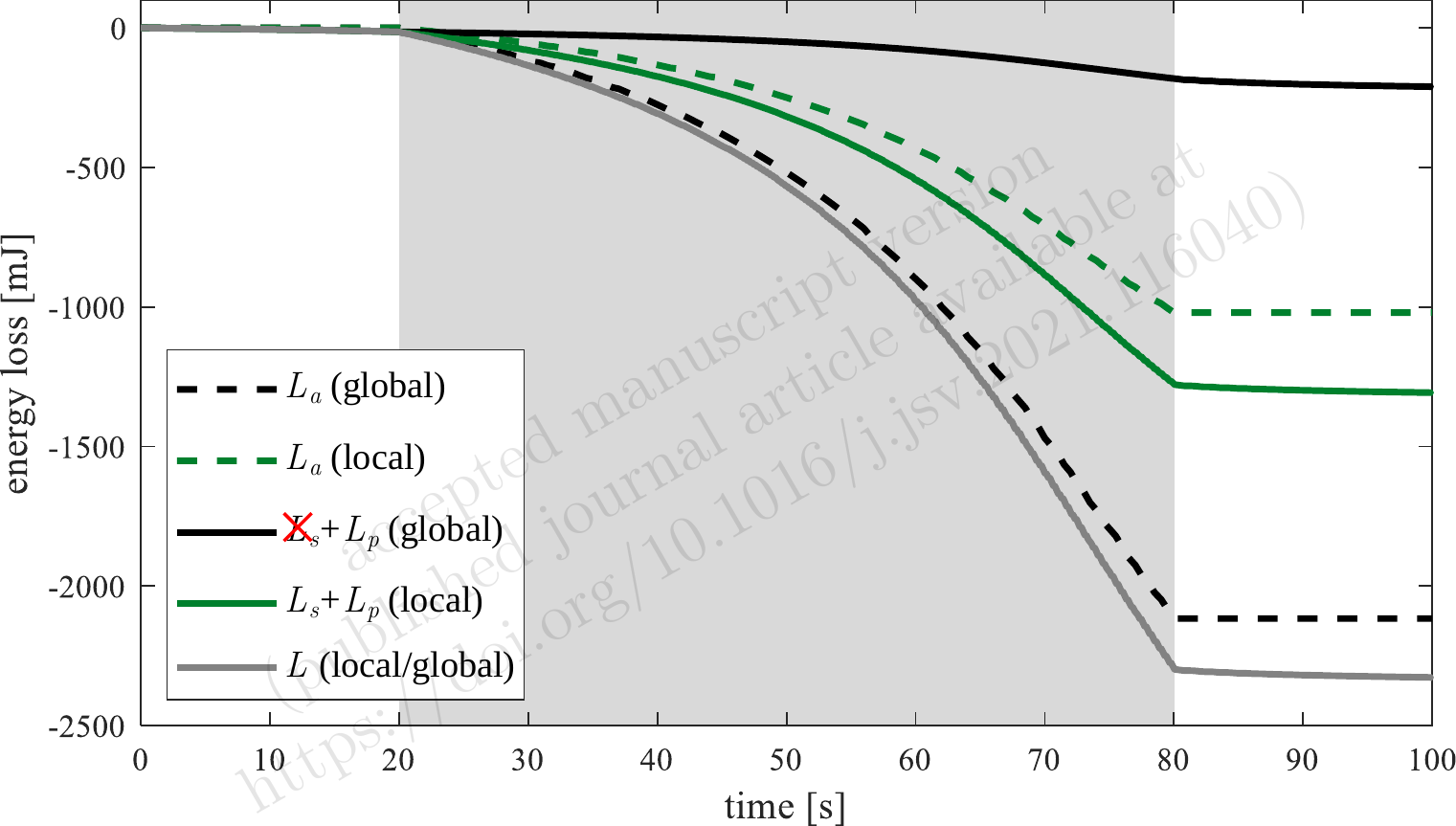}
	\caption{Energy loss, separated into pseudo-active and semi-active + passive component (sweep)}
	\label{fig:11}
\end{figure}

With this modal consideration, the individual contributions to the energy loss ($ L_a $, $ L_s $, and $ L_p $, see \autoref{subsec:2.4}) can now be accounted for.
The above analyzed dissipative energy loss $ L_s + L_p $ \autoref{eq:39} is now considered for the whole system (mode~1 and mode~2), along with the extracted energy $ L_a $ \autoref{eq:38} for the local and global stiffness variation (see \autoref{fig:11}).
As mentioned, the total energy loss $ L $ \autoref{eq:37} (gray curve) is identical for both cases.
In the global case, the energy loss consists mainly of extracted energy (pseudo-active effect, dashed black curve).
The small remaining portion corresponds to the passive part $ L_p $ (passive effect, solid black curve).
This portion results exclusively from the dissipated energy of mode~1 (see \autoref{fig:10}, dashed black curve).
Thus, no semi-active component $ L_s $ is present, and the dissipated energy $ L_s + L_p $ is reduced to $ L_p $.
In the local case, the energy transfer to the second mode results in a substantial semi-active contribution, which, together with the passive one (both combined in the solid green curve), exceeds the pseudo-active component (dashed green curve) for the selected parameter setting.

In summary, it can be stated that for harmonic excitation, resonances can be effectively prevented using the presented method of stiffness modulation.
If stiffness adaptation is possible, however, vibrations can also be effectively reduced by a constant stiffness change which moves the system's eigenfrequency away from the excitation frequency (see \autoref{subsec:1.1});
however, this would only work under ideal conditions, such as a known narrowband excitation.
With the cyclic variations of the structure's stiffness based on the amplitude of the considered mode, we are not restricted to such simplified situations.
Furthermore, if local stiffness modulation is applied, this leads to an internal energy transfer and thus to a higher ratio of energy loss due to inherent damping to energy extraction by the actuator.

\section{Conclusions and outlook}
\label{sec:3}
In this paper, a study on semi-active vibration reduction on the basis of cyclic stiffness variations (modulation) was presented, with focus on the so-called switched-stiffness approach.
Compared with previous work on this topic, the novelty of this contribution consists in the discrimination of the physical mechanism of vibration reduction into a purely semi-active effect, which enhances dissipation within the structure by redistributing energy among modes, and a pseudo-active effect, which extracts vibration energy from the system by negative work of the stiffness variation device.
Both contributions join the passive effect, i.e., the energy dissipation that would be present without stiffness modulation.

The semi-active effect is expected to be more advantageous because it does not require mechanical work to be performed by the stiffness variation device, which reduces the hardware size and cost.
This energetically passive nature also fulfills the common assumption formulated for semi-active measures.
To exploit it to the fullest, it is crucial to deepen how the semi-active share of vibration attenuation can be influenced by design.
The analysis of spring-mass oscillators showed that the spatial distribution of the stiffness variation plays a central role in this sense.
The semi-active effect is only present when the stiffness is changed locally and disappears when the stiffness change is homogeneous over the system.
The exploitation of the semi-active effect allows for the more effective use of the inherent damping capacity of structures, and reduces the need for cumbersome external devices.
This will be of particular interest for lightweight systems.

In future numerical and experimental studies, concrete options for stiffness variation devices will have to be investigated.
Variable-stiffness beam structures can be realized, for instance, by shape adaption of the cross section, which has a strong effect on the bending and torsion stiffness.
A thin-walled construction with several independently controllable shape-adaptable ribs offers the possibility of differentiated actuation.
In this way, the positive effect of local stiffness variation on the semi-active effect can be verified for continuous systems.
Another option for the physical realization of stiffness variation consists in the use of bending elements with variable prestress.

Other topics for future studies are provided by alternative choices of the control logic.
In the presence of a physical stiffness variation device, the steps included in the stiffness-switching strategy involve impulsive excitation of the systems, which can lead to unwanted vibrations.
For this reason, it can be appropriate to move from the theoretical optimum control logic to a smoother time law.
Further, continuous systems offer a theoretically infinite range of choices for the observation function; thus, dedicated studies are needed to identify the best options.
Finally, the above-mentioned extensions will be considered for additional types of excitation of the vibration structure.

\section*{Acknowledgements}
Funded by the Deutsche Forschungsgemeinschaft (DFG, German Research Foundation) -- SPP~1897 `Calm, Smooth and Smart' (project numbers HA 7893/1-1 and WI 1181/10-1).

\printbibliography
\end{document}
\endinput